\documentclass[11pt,a4paper]{article}
\usepackage[T1]{fontenc}
\usepackage{lmodern}
\usepackage[margin=26mm,headheight=14pt]{geometry}
\usepackage{amsmath,amssymb,amsthm,mathtools,bm}
\usepackage{microtype,booktabs,array,graphicx,needspace,flafter,placeins}
\usepackage{xcolor}
\usepackage[hidelinks]{hyperref}
\usepackage{fancyhdr}
\allowdisplaybreaks[1]
\numberwithin{equation}{section}
\newtheorem{theorem}{Theorem}[section]
\newtheorem{lemma}[theorem]{Lemma}
\newtheorem{proposition}[theorem]{Proposition}
\newtheorem{corollary}[theorem]{Corollary}
\theoremstyle{definition}
\newtheorem{assumption}[theorem]{Assumption}

\theoremstyle{remark}
\newtheorem{remark}[theorem]{Remark}
\DeclareMathOperator{\diver}{div}
\DeclareMathOperator{\spec}{spec}
\DeclareMathOperator{\rad}{rad}
\DeclareMathOperator{\dist}{dist}
\DeclareMathOperator{\Rea}{Re}
\DeclareMathOperator{\Ima}{Im}
\newcommand{\dd}{\,\mathrm d}
\newcommand{\Om}{\Omega}

\newcommand{\Tcal}{\mathcal T}
\newcommand{\Dcal}{\mathcal D}
\newcommand{\Kcal}{\mathcal K}
\newcommand{\Gcal}{\mathcal G}
\newcommand{\Ucal}{\mathcal U}
\newcommand{\Pcal}{\mathcal P}
\newcommand{\norm}[1]{\left\lVert#1\right\rVert}

\newcommand{\wpi}{\varpi}

\newcommand{\bR}{\mathbb R}

\newcommand{\qe}{q_{\mathrm e}}

\newcommand{\epsN}{\varepsilon_N}
\newcommand{\invN}{C_{\mathrm{inv},N}}
\newcommand{\mean}[1]{\langle#1\rangle_{\Om}}

\AddToHook{env/theorem/before}{\Needspace{4\baselineskip}}
\AddToHook{env/lemma/before}{\Needspace{3\baselineskip}}
\AddToHook{env/proposition/before}{\Needspace{3\baselineskip}}
\AddToHook{env/corollary/before}{\Needspace{3\baselineskip}}
\usepackage{algorithm,algpseudocode}
\AddToHook{env/assumption/before}{\Needspace{4\baselineskip}}
\newcommand{\InputTableRows}[1]{\csname @@input\endcsname #1 }
\hypersetup{pdftitle={Flux-Matched Poisson Iteration for 
Generated Jacobian Equations: 
Local Convergence and Spectral Galerkin Approximation}}
\title{\vspace{-8mm}Flux-Matched Poisson Iteration 
for\\Generated Jacobian Equations:\\[2mm] 
Local Convergence and Spectral Galerkin Approximation}
\newcommand{\AuthorOne}{Haoxiang Zhang}
\newcommand{\AuthorTwo}{Yixuan Zhang}
\newcommand{\AffiliationOne}{School of Mathematical Sciences, Zhejiang University, Hangzhou, 310027, China}
\newcommand{\AffiliationTwo}{School of Mathematical Science, Peking University, Beijing, 100871, China}
\newcommand{\EmailOne}{zhanghx24@zju.edu.cn}
\newcommand{\EmailTwo}{yixuan@math.pku.edu.cn, yixuan@pku.edu.cn}

\author{%
  \AuthorOne\textsuperscript{1},\quad
  \AuthorTwo\textsuperscript{2}\\[3pt]
  {\small\textsuperscript{1}\AffiliationOne}\\[1pt]
  {\small\textsuperscript{2}\AffiliationTwo}\\[2pt]
  {\small Emails: \texttt{\EmailOne};\quad
    \texttt{\EmailTwo}}}
\hypersetup{pdfauthor={\AuthorOne; \AuthorTwo}}

\date{}
\begin{document}
\maketitle
\vspace{-8mm}
\begin{abstract}
We develop a flux-matched Poisson iteration for solving
generated Jacobian equations with second boundary conditions.
In each step, it solves a fixed shifted Poisson problem with a Neumann-type
boundary condition. 
Near a smooth regular solution with positive $C^{1,1}$ densities, 
we prove that the iteration has local geometric convergence 
in $C^{2,\alpha}$ on smooth connected Riemannian domains. 
The iteration is then discretized using a spectral Galerkin method, 
with local convergence proved for the Fourier--Jacobi scheme on spherical caps
Optical numerical experiments support these estimates and demonstrate 
useful imaging performance beyond the local convergence theory.
\end{abstract}

\section{Introduction}
\label{sec:introduction}\label{sec:opticalmodels}
Generated Jacobian equations (GJEs) extend the Monge--Amp\`ere equation through 
maps generated by an unknown potential. Let $\Omega,\Lambda$ be smooth Riemannian source and target domains 
of the same dimension at least two, with volume measures $\dd V,\dd V_\Lambda$. 
The scalar potential $u:\Omega\to\mathbb R$ satisfies
\begin{equation}
 \det\!\left[\nabla^2u-A(x,u,\nabla u)\right]
 =B(x,u,\nabla u),
 \label{eq:generalGJE}
\end{equation}
where the symmetric $A$ and scalar $B$ arise from a generating function 
and prescribed data ~\cite{Trudinger}.

A smooth generator $\Gcal(x,y,z)$, with source point $x$, target point $y$, and auxiliary scalar $z$, 
determines $(T_u,z_u)$ through
\begin{equation}
 u(x)=\Gcal(x,T_u(x),z_u(x)),\quad
 \nabla u(x)=\Gcal_x(x,T_u(x),z_u(x)).
 \label{eq:contact}
\end{equation}
The equations are locally solvable and define $T_u(x)=Y(x,u(x),\nabla u(x))$. 
Evaluating generator derivatives at contact, set
\begin{equation}
 A_u:=A(x,u,\nabla u)=\Gcal_{xx},\quad
 E_u:=\Gcal_{xy}
       -\frac{\Gcal_{xz}\otimes\Gcal_y}{\Gcal_z}.
 \label{eq:coeffs}
\end{equation}
The metrics identify $E_u$ as a target-to-source tangent map. 
Conditions $\Gcal_z\ne0$ and $\det E_u\ne0$ ensure local solvability. 
Differentiating $u$ twice with (\ref{eq:coeffs}) gives $E_uDT_u=W_u:=\nabla^2u-A_u$.

The positive volume Jacobian satisfies $T_u^*(\dd V_\Lambda)=\lvert\det(\nabla T_u)\rvert\dd V$. 
On the elliptic branch $W_u\succ0$, taking determinants in $E_uDT_u=W_u$ gives
\begin{equation}
 \lvert\det(\nabla T_u)\rvert=\frac{\det(\nabla^2u-A_u)}{|\det E_u|}.
 \label{eq:areaJacobian}
\end{equation}
The absolute value allows either orientation of the map.

We consider the GJEs in the mass-transport class. 
For positive $C^{1,1}$ densities $f$ on $\overline\Omega$ and $g$ on $\overline\Lambda$ of equal mass,
seek $u$ transporting $f\dd V$ onto $g\dd V_\Lambda$. 
Local conservation $g(T_u)\lvert\det(\nabla T_u)\rvert=f$ with (\ref{eq:areaJacobian}) gives
\begin{equation}
 \det(\nabla^2u-A_u)
 =|\det E_u|\frac{f(x)}{g(T_u(x))},\quad W_u\succ0.
 \label{eq:generatedMA}
\end{equation}
The second boundary condition prescribes the target image. 
For $T_u$ a diffeomorphism of closures, problem (\ref{eq:generatedMA}) is equivalent to
\begin{equation}
 g(T_u)\lvert\det(\nabla T_u)\rvert=f\quad\text{in }\Omega,\quad
 T_u(\partial\Omega)=\partial\Lambda,\quad T_u(\Omega)=\Lambda.
 \label{eq:physicalproblem}
\end{equation}
We express this problem through its mass and boundary residuals.
Extend $g$ positively to a $C^{1,1}$ function $g_e$ near $\overline\Lambda$,
and let $\psi$ be signed boundary distance, negative inside.
For $0<\alpha<1$ and a fixed positive weight $a\in C^{1,\alpha}(\partial\Omega)$, define
\begin{equation}
 F(u)=f-g_e(T_u)\lvert\det(\nabla T_u)\rvert,\quad
 \Psi(u)=a\,\psi(T_u)|_{\partial\Omega}.
 \label{eq:residuals}
\end{equation}
Solutions of \eqref{eq:physicalproblem} satisfy $F=\Psi=0$.
Classical optimal transport uses $\Gcal=-c(x,y)-z$, where $c$ is the transport cost.
In this case, $T_u$ is independent of the potential value. 
General generators also introduce dependence on $u$, 
producing drift in the linearized equation. 
In optics, $T_u$ describes outgoing rays or receiver intersections.

GJE regularity rests on structural and domain assumptions. 
Trudinger~\cite{Trudinger,TrudingerRevisited} developed the local theory. 
Guillen and Kitagawa~\cite{GuillenKitagawa} proved strict generated 
convexity and interior regularity. 
Second boundary value solvability was established for cost-generated 
equations~\cite{TrudingerWang} and general generators~\cite{JiangTrudinger}, 
with further global regularity and derivative estimates in~\cite{RankinGlobal,TrudingerSecond}. 
We assume a smooth regular solution and analyze a nearby iteration.

For finite targets, convergent GJE algorithms include coordinatewise mass 
adjustment~\cite{AbedinGutierrez} and damped Newton iteration~\cite{GallouetMerigotThibert}. 
PDE methods include GJE least squares~\cite{RomijnEtAl}, 
reflector collocation~\cite{BrixHafizogullariPlaten}, 
and convergent finite differences for spherical transport and 
reflection~\cite{FroeseSphere,FroeseHamfeldtTurnquist}. 
Poisson-based Monge--Amp\`ere solvers include time marching~\cite{Awanou}, 
a recent $L$-scheme~\cite{KohleEtAl}, and least-squares splitting~\cite{PerusoSorella}. 
Spectral Galerkin methods are studied in~\cite{JinEtAl}. Zhang and Bao~\cite{ZhangBao} 
apply Sobolev gradient descent to far-field reflector design.

Inspired by this far-field Sobolev gradient descent approach, we propose a
flux-matched Poisson iteration for solving GJEs with second boundary conditions.
Each step solves a fixed shifted Poisson problem driven by the mass and boundary residuals. 
A reference-matched boundary weight combines their linearizations into a conservative 
diffusion--drift flux. Under the stated regularity and parameter assumptions, 
we prove local geometric convergence of the fixed-point iteration in $C^{2,\alpha}$ and decay of
both residuals near a smooth regular solution. The drift-free case permits zero shift.

On spherical caps, we construct a Fourier--Jacobi Galerkin scheme whose fixed matrix 
can be factored once. With exact integration, 
we prove local convergence to nearby discrete fixed points at a rate uniform 
in degree within possibly shrinking neighborhoods. 
The estimates separate iteration and spatial errors and 
yield superalgebraic spatial convergence for smooth references. 
Numerical tests examine smooth convergence and the boundary correction. 
Reflector imaging and square-domain refraction illustrate practical performance 
beyond the proved setting.

Section~\ref{sec:algorithm} defines the iteration and states the continuous convergence theorem. 
Section~\ref{sec:continuous} gives its proof. 
Section~\ref{sec:spectral} constructs the spectral discretization. 
Section~\ref{sec:discrete} proves discrete convergence and error bounds. 
Section~\ref{sec:numerics} presents the numerical experiments. 
Section~\ref{sec:conclusion} concludes the paper.

\label{page:introduction-end}
\section{The Poisson iteration and local convergence}
\label{sec:geometry}\label{sec:algorithm}

We construct a flux-matched Poisson iteration for the residual problem
\eqref{eq:residuals} and state its local convergence theorem.

\subsection{Flux-matched Poisson correction}
\label{sec:iteration}
Given a zero-mean $C^{2,\alpha}$ GJE approximation $u^n$, we compute the
correction $v^n$ using a fixed shifted Poisson solver.
The mass residual $F(u^n)$ provides the interior source, and the weighted
boundary residual $\Psi(u^n)$ provides the Neumann data.
The next iterate is $u^{n+1}=u^n-\tau v^n$.

For fixed parameters $\sigma\ge0$ and $\tau>0$,
the iteration is
\begin{equation}
 \begin{cases}
  (\sigma I-\Delta)v^n=F(u^n)-m(u^n)&\text{in }\Omega,\\
  \partial_\nu v^n=\Psi(u^n)&\text{on }\partial\Omega,\\
  \displaystyle\int_\Omega v^n\dd V=0,
 \end{cases}
 \quad u^{n+1}=u^n-\tau v^n.
 \label{eq:continuousmap}
\end{equation}
Here $\nu$ is the outward conormal. 
Integration with the boundary data and zero mean gives
\begin{equation}
 m(u)=\frac1{|\Omega|}
 \left(\int_\Omega F(u)\dd V+
       \int_{\partial\Omega}\Psi(u)\dd S\right).
 \label{eq:meanresidual}
\end{equation}
It enforces Neumann compatibility
at $\sigma=0$ and zero mean at $\sigma>0$. 
The update can be denoted as a map $S(u^n) = u^{n+1}$. 
The fixed weight in $\Psi$ converts boundary displacement into normal flux.

\Needspace{4\baselineskip}
\subsection{Regular reference}
Each correction has zero mean, so the iteration preserves the initial mean.
We use zero-mean variables for the analysis. To center a physical potential,
we subtract its mean from both $u$ and $\Gcal$. This preserves the generated
map even when it depends on the potential value.

\begin{assumption}[Regular reference]\label{ass:reference}
\leavevmode
\begin{enumerate}
\item We measure local convergence of the iteration method in 
the $C^{2,\alpha}$ norm. 
For $0<\alpha<1$, define the zero-mean state space
\begin{equation}
 X=\left\{h\in C^{2,\alpha}(\overline\Omega):
                   \int_\Omega h\dd V=0\right\},
 \quad \norm h_X=\norm h_{C^{2,\alpha}(\overline\Omega)}.
 \label{eq:spaces}
\end{equation}
\item The closures $\overline\Omega$ and $\overline\Lambda$ are connected,
compact, smooth Riemannian manifolds of the same dimension at least two,
with smooth boundaries. The target lies in a smooth ambient manifold
supporting nearby maps and boundary collars. The strictly positive densities
satisfy $f\in C^{1,1}(\overline\Omega)$ and $g\in C^{1,1}(\overline\Lambda)$. Assume
\begin{equation}
 \int_\Omega f\dd V=\int_\Lambda g\dd V_\Lambda.
 \label{eq:eqal mass}
\end{equation}
\item There exists
$u^\ast\in C^\infty(\overline\Omega)\cap X$ satisfying
\eqref{eq:physicalproblem}, such that
$T_*=T_{u^\ast}:\overline\Omega\to\overline\Lambda$ is a diffeomorphism.
The generator is smooth near the compact reference contact graph, where
$\Gcal_z\ne0$, $E_*:= E_{u^*}$ is invertible, 
and $W_* := W_{u^*}\succ0$ uniformly. 
The reference tensor $K=fW_*^{-1}$ is real and belongs to $C^{1,1}$, with
$0<\kappa_-\le\kappa_+$ such that for 
every $x \in \overline{\Omega}$, 
\begin{equation}
 \kappa_-I\preceq K(x)=K(x)^*\preceq\kappa_+I.
 \label{eq:ellipticity}
\end{equation}
\end{enumerate}
\end{assumption}

In addition to the diffusion tensor $K$, the reference solution determines
the vector field
\begin{equation}
 d:=d_{u^\ast},\quad
 d_u(x):=\frac{\Gcal_{xz}}{\Gcal_z}\bigl(x,T_u(x),z_u(x)\bigr).
 \label{eq:referencedrift}
\end{equation}
As shown in Section~\ref{sec:matching}, this coefficient enters the
linearized flux as $K(\nabla h-dh)$. We call the reference drift-free when
$d\equiv0$. 

On $\partial\Omega$, define the reference-matched weight by
\begin{equation}
 n_*=(\nabla\psi)\circ T_*,\quad
 \chi=n_*\cdot DT_*\nu>0,\quad a_*=f/\chi.
 \label{eq:referenceweight}
\end{equation}
Boundary correspondence and the interior property give $\chi>0$.

\subsection{The continuous local convergence theorem}
\begin{theorem}[Continuous local convergence]\label{thm:continuous}
Under Assumption~\ref{ass:reference}, use the matched weight $a_*$ from
\eqref{eq:referenceweight} for $\Psi$ in \eqref{eq:continuousmap}.
There are finite constants $\sigma_*\ge0$ and
$\tau_*>0$, determined by the reference coefficients and source geometry,
such that for every fixed 
\begin{equation}
 \sigma\ge\sigma_*,\quad 0<\tau<\tau_*,
 \label{eq:parameters}
\end{equation}
there exist $\delta>0$, $C\ge1$, and $q\in(0,1)$ with the following property.
Every $u^0\in X$ satisfying $\norm{u^0-u^\ast}_X<\delta$ generates a
well-defined trajectory of (\ref{eq:continuousmap}) 
that stays in a regular neighborhood of $u^\ast$ and obeys
\begin{equation}
 \norm{u^n-u^\ast}_X\le Cq^n\norm{u^0-u^\ast}_X.
 \label{eq:continuousconvergence}
\end{equation}
The mass, boundary, and compatibility residuals satisfy
\begin{equation}
 \norm{F(u^n)}_{C^{0,\alpha}}+
 \norm{\Psi(u^n)}_{C^{1,\alpha}(\partial\Omega)}+|m(u^n)|
 \le Cq^n\norm{u^0-u^\ast}_X.
 \label{eq:continuousresiduals}
\end{equation}
If the reference drift vanishes, $d=0$, the admissible range can instead be
taken as
\begin{equation}
 \sigma\ge0,\quad 0<\tau<2/\kappa_+.
 \label{eq:driftfree}
\end{equation}
\end{theorem}

Constants and neighborhoods may depend on the dimension and fixed parameters.
The transport boundary condition is recovered in the limit 
rather than imposed at every step.

\section{Proof of continuous local convergence}
\label{sec:continuous}\label{sec:main}\label{sec:linear}

To prove the main theorem, 
we first identify the matched linearized flux and use a positive stationary
density to estimate its eigenvalues. Compact perturbation and conormal
resolvent arguments then control the full spectrum on $X$.
An equivalent norm yields a locally invariant contraction neighborhood.

We begin with the energy space of the linearized iteration.
Testing \eqref{eq:continuousmap} with a zero-mean function removes $m(u^n)$.
Integrating by parts leads to
\begin{equation}
 \begin{aligned}
 H&=\left\{h\in H^1(\Omega):\int_\Omega h\dd V=0\right\},\\
 p_\sigma(w,\phi)&=\int_\Omega\nabla w\cdot\nabla\overline\phi\dd V
                  +\sigma\int_\Omega w\overline\phi\dd V,
 \quad \norm h_\sigma^2=p_\sigma(h,h).
 \end{aligned}
 \label{eq:energy}
\end{equation}
The weak correction equation of (\ref{eq:continuousmap}) is
\begin{equation}
 p_\sigma(v^n,\phi)=\int_\Omega F(u^n)\overline\phi\dd V
       +\int_{\partial\Omega}\Psi(u^n)\overline\phi\dd S,
 \quad \phi\in H.
 \label{eq:iterationweak}
\end{equation}
Poincar\'e's inequality makes $\norm\cdot_\sigma$ equivalent to the $H^1$ norm,
including at $\sigma=0$. The correction $v^n$ is the Riesz representative
of the combined residual load in the $p_\sigma$ inner product.

\Needspace{4\baselineskip}
\subsection{Linearization and the matched flux}
\label{sec:matching}
\begin{lemma}[Neighborhood regularity]\label{lem:neighborhood}
Under Assumption~\ref{ass:reference}, there is an open neighborhood $\Ucal\subset X$
on which the chosen contact branch is well defined. 
For every $u\in\mathcal U$, the generated map $T_u$
is an embedding into the ambient target manifold, and $W_u\succ0$.
The residuals are well defined there. In local coordinates, the maps 
\begin{equation}
 u\mapsto T_u:\Ucal\to C^{1,\alpha}(\overline\Omega),\quad
 F:\Ucal\to C^{0,\alpha}(\overline\Omega),\quad
 \Psi:\Ucal\to C^{1,\alpha}(\partial\Omega)
 \label{eq:regularity}
\end{equation}
are continuously differentiable.
\end{lemma}
\begin{proof}
The contact Jacobian $\Gcal_z\det E\ne0$ gives a smooth branch
$T_u=Y(x,u,\nabla u)$ by the implicit function theorem \cite{Trudinger}.
Compactness makes the neighborhood uniform. H\"older product estimates
give smooth dependence of $T_u$, its nonvanishing Jacobian, and $\Psi$
in the stated spaces. Since $\nabla g_e$ is Lipschitz and $T_u$ is
uniformly Lipschitz in space, H\"older interpolation gives continuity of
$\nabla g_e(T_u)$ in $C^{0,\alpha}$ with modulus
$O(\norm{u-\widetilde u}_X^{1-\alpha})$. Integrating the first variation
of $F$ gives a remainder $O(\norm h_X^{2-\alpha})$, so $F$ is $C^1$.
Openness of embeddings in $C^1$ \cite[Chapter~2]{Hirsch} and continuity
preserve the remaining conditions after shrinking $\Ucal$ to a convex
neighborhood.
\end{proof}

\begin{lemma}[Neumann correction]\label{lem:solver}
Fix $\sigma\ge0$. For an interior source $F_0\in C^{0,\alpha}(\overline\Omega)$
and Neumann data $b_0\in C^{1,\alpha}(\partial\Omega)$, set
$c_0=|\Omega|^{-1}(\int_\Omega F_0\dd V+\int_{\partial\Omega}b_0\dd S)$.
The problem
\begin{equation}
 (\sigma I-\Delta)v=F_0-c_0\quad\text{in }\Omega,\quad
 \partial_\nu v=b_0\quad\text{on }\partial\Omega,\quad
 \int_\Omega v\dd V=0
 \label{eq:neumannproblem}
\end{equation}
has a unique solution $v\in X$. It satisfies
\begin{equation}
 p_\sigma(v,\phi)=\int_\Omega F_0\overline\phi\dd V+
                  \int_{\partial\Omega}b_0\overline\phi\dd S,
 \quad \phi\in H,
 \label{eq:weaksolver}
\end{equation}
and $\norm v_X\le C_\sigma(\norm{F_0}_{C^{0,\alpha}}+
\norm{b_0}_{C^{1,\alpha}})$.

For fixed $a$ and $\tau>0$, take $F_0=F(u)$ and $b_0=\Psi(u)$, and denote
the corresponding solution by $v=\varphi(u)$. Then $c_0=m(u)$, and
$\varphi:\Ucal\to X$, $S=I-\tau\varphi:\Ucal\to X$, and
$m:\Ucal\to\mathbb R$ are $C^1$. On a smaller regular neighborhood, 
$D\varphi$ is locally Lipschitz in $\mathcal{L}(X, H)$.
\end{lemma}
\begin{proof}
Poincar\'e's inequality and the trace theorem allow Lax--Milgram to give
a unique weak solution in $H$, with an $H^1$ bound. Testing with
$\phi-\mean\phi$ for arbitrary $\phi\in H^1$ identifies
\eqref{eq:neumannproblem}. Boundary regularity for this Neumann problem
gives $v\in X$. The Schauder estimate, with its lower-order term controlled
by the $H^1$ bound, proves the stated estimate \cite{ADN1,ADN2,Nardi}.

The solution depends linearly and continuously on the data, so
Lemma~\ref{lem:neighborhood} gives the $C^1$ claims. Since $\nabla g_e$
is Lipschitz and the geometric factors depend smoothly on $u$,
$DF$ and $D\Psi$ are locally Lipschitz as operators from $X$ into the
interior and boundary $L^2$ spaces. Subtract the differentiated weak
equations and test with $z=(D\varphi(u)-D\varphi(\widetilde u))h$.
The energy and trace estimates prove local Lipschitz continuity, also
at $\sigma=0$. 
\end{proof}

We next linearize the generated map.
Write $\delta T=D_uT(u)[h]$,
where $DT_u$ denotes the spatial differential of the generated map.
Varying and differentiating \eqref{eq:contact} gives
$E_u\delta T=\nabla h-d_uh$, with $d_u$ defined in
\eqref{eq:referencedrift}, and $E_uDT_u=W_u$. Hence
\begin{equation}
 (DT_u)^{-1}\delta T=W_u^{-1}(\nabla h-d_uh).
 \label{eq:mapdiff}
\end{equation}
For $\rho_u=g_e(T_u)\lvert\det(\nabla T_u)\rvert$, the Piola variation identity is
$D\rho_u[h]=\diver\bigl(\rho_u(DT_u)^{-1}\delta T\bigr)$.
Using $F(u) = f - \rho_u$ and (\ref{eq:mapdiff}), we obtain
\begin{equation}
 DF(u)h=-\diver\bigl(\rho_uW_u^{-1}(\nabla h-d_uh)\bigr).
 \label{eq:DF}
\end{equation}
For $u\in C^{2,\alpha}$, this identity is understood distributionally. 

At the reference, $\rho_{u^\ast}=f$. Define the flux and operator
\begin{equation}
 j_h=K(\nabla h-dh),\quad \mathcal Lh=-\diver j_h.
 \label{eq:linearizedflux}
\end{equation}

\begin{proposition}[Boundary--flux matching]\label{prop:matching}
For $\Psi(u)=a_*\psi(T_u)|_{\partial\Omega}$ with $a_*$ defined in
\eqref{eq:referenceweight},
\begin{equation}
 DF(u^\ast)h=-\diver j_h,\quad
 D\Psi(u^\ast)h=\nu\cdot j_h.
 \label{eq:pairderivative}
\end{equation}
\end{proposition}
\begin{proof} 
$DF(u^*)$ follows from \eqref{eq:DF} and \eqref{eq:linearizedflux}. 
Boundary correspondence gives $(DT_*)^*n_*=\chi\nu$ for the quantities
defined in \eqref{eq:referenceweight}.
Combining it with \eqref{eq:mapdiff} at $u^\ast$ and
$a_*\chi=f$, we obtain
\[
 D\Psi(u^\ast)h=a_*n_*\cdot\delta T
 =a_*\chi\nu\cdot W_*^{-1}(\nabla h-dh)
 =\nu\cdot j_h.
\]
\end{proof}

The divergence theorem gives
$Dm(u^\ast)h=0$, so
$A_\sigma=D\varphi(u^\ast)$ has the weak form
\begin{equation}
 p_\sigma(A_\sigma h,\phi)
 =\int_\Omega K(\nabla h-dh)\cdot\nabla\overline\phi\dd V,
 \quad \phi\in H.
 \label{eq:weakA}
\end{equation}
It defines a bounded extension of $A_\sigma$ to $H$. On space $X$, the flux
in \eqref{eq:linearizedflux} gives the strong form
\begin{equation}
 (\sigma I-\Delta)A_\sigma h=-\diver j_h,
 \quad\partial_\nu A_\sigma h=\nu\cdot j_h,
 \quad\int_\Omega A_\sigma h\dd V=0.
 \label{eq:strongA}
\end{equation}
The linearized update is therefore
\begin{equation}
 \Tcal=DS(u^\ast)=I-\tau A_\sigma.
 \label{eq:updateDerivative}
\end{equation}

\subsection{Eigenvalue estimates}

We next estimate the eigenvalues of $A_\sigma$. The weak form of $A_\sigma$
contains a symmetric diffusion term and a drift term whose
contribution to the energy has no definite sign. To handle the
drift, we construct a positive stationary density $\wpi$.
Set $h=\wpi w$. It eliminates the stationary flux contribution
from the real part of the weighted energy identity.

\begin{lemma}[Positive reflecting density]\label{lem:density}
For real $K,d\in C^{1,\alpha}$ satisfying \eqref{eq:ellipticity}, the problem
\begin{equation}
 \mathcal L\wpi=0 \text{ in }\Omega,\quad
 \nu\cdot K(\nabla\wpi-d\wpi)=0 \text{ on }\partial\Omega,
 \quad\int_\Omega\wpi\dd V=1
 \label{eq:invariant}
\end{equation}
has a unique solution $\wpi\in C^{2,\alpha}(\overline\Omega)$ that is strictly
positive on the closure. The weak kernel of this no-flux problem on the
full $H^1$ space is $\operatorname{span}\{\wpi\}$.
\end{lemma}
\begin{proof}
Choose $\eta>\|Kd\|_\infty^2/(2\kappa_-)$ sufficiently large that
$\eta+\diver(Kd)>0$, and consider the form on $H^1(\Omega)$ 
\[
 \mathfrak a_\eta(\zeta,\phi)=\int K\nabla\zeta\cdot\nabla\overline\phi
 -\int \zeta Kd\cdot\nabla\overline\phi+\eta\int\zeta\overline\phi.
\]
Young's inequality gives
\begin{equation}
 \Rea\mathfrak a_\eta(\zeta,\zeta)\ge\frac{\kappa_-}{2}\|\nabla\zeta\|_2^2+
 \left(\eta-\frac{\|Kd\|_\infty^2}{2\kappa_-}\right)\|\zeta\|_2^2.
 \label{eq:densitycoercivity}
\end{equation}
By Lax-Milgram and Rellich, the inverse map $R_\eta=(\eta+\mathcal L)^{-1}:L^2\to H^1$ 
with the natural no-flux condition is well-defined and compact on $L^2$. 

For $r\ge0$ and $\zeta=R_\eta r$, testing with $\zeta_-=\max\{-\zeta,0\}$
and using $\nabla \zeta_-=-\mathbf1_{\{\zeta<0\}}\nabla \zeta$ give 
\[
 0\le\int r\zeta_-=\mathfrak a_\eta(\zeta,\zeta_-)
                  =-\mathfrak a_\eta(\zeta_-,\zeta_-).
\]
Coercivity \eqref{eq:densitycoercivity} forces $\zeta_-=0$, so $R_\eta$ is positive. 
We next construct a positive eigenfunction of $R_\eta$ 
corresponding to the eigenvalue $\eta^{-1}$ and prove it is the density we want.

Integration by parts gives the adjoint operator 
$\mathcal L^*\phi=-\diver(K\nabla\phi)-Kd\cdot\nabla\phi$ and $\nu\cdot K\nabla\phi=0$. 
$R_\eta^*$ is also positive by duality, and
$\mathcal L^*1=0$ gives $R_\eta^*1=\eta^{-1}1$. Thus
\[
 |R_\eta^*\phi|\le R_\eta^*|\phi|
 \le\|\phi\|_\infty R_\eta^*1
 =\eta^{-1}\|\phi\|_\infty.
\]
For an $L^2$ eigenvector $R_\eta^*\phi=\mu\phi$ with $\mu\ne0$, it holds 
$(\eta+\mathcal L^*)\phi=\mu^{-1}\phi$ with the adjoint conormal condition. 
Repeated application of boundary elliptic estimates and
Sobolev embedding yield 
$\phi\in C(\overline\Omega)$, 
so the displayed estimate gives $|\mu|\le\eta^{-1}$.
Since $R_\eta^*$ is compact, this bounds every nonzero
spectral value. The constant eigenfunction attains the bound, 
hence $r_{L^2}(R_\eta)=r_{L^2}(R_\eta^*)=\eta^{-1}>0$.
The positive compact-operator theorem~\cite{Schaefer} gives a nonzero $\wpi\ge0$
with $R_\eta\wpi=\eta^{-1}\wpi$, hence $\mathcal L\wpi=0$ with no flux.

The condition $\nu\cdot K\nabla\wpi=(\nu\cdot Kd)\wpi$
is strictly oblique since $\nu\cdot K\nu\ge\kappa_-$.
Boundary $L^p$ bootstrapping
followed by Schauder estimates gives
$\wpi\in C^{2,\alpha}(\overline\Omega)$ \cite{ADN1,ADN2}.
For strict positivity, consider the operator
\[
 P_\eta=-\mathcal L-\eta
 =\diver(K\nabla\cdot)-Kd\cdot\nabla-(\eta+\diver(Kd)),
\]
which has a nonpositive zeroth-order coefficient and satisfies
$P_\eta\wpi=-\eta\wpi\le0$.
The strong minimum principle gives $\wpi > 0$ in $\Omega$.
At a boundary zero $x_0$, tangential derivatives vanish and the no-flux
condition gives $(\nu\cdot K\nu)\partial_\nu\wpi=\nu\cdot K\nabla\wpi=0$.
It implies $\partial_\nu\wpi(x_0)=0$ and contradicts the Hopf lemma \cite{GT}, 
which gives $\partial_\nu\wpi(x_0)<0$. 
Thus $\wpi>0$ on $\overline\Omega$ and can be normalized to integral one. 

For any weak nullvector $h\in H^1$, set $w=h/\wpi$ and
$J=K(\nabla\wpi-d\wpi)$. Since $\min\wpi>0$, it holds $w\in H^1$.
Using the flux definition \eqref{eq:linearizedflux}, we obtain
\begin{align}
 j_h=\wpi K\nabla w+wJ,\quad
 &\diver J=0,\quad J\cdot\nu=0,\nonumber\\
 \Rea\int wJ\cdot\nabla\overline w
 &=\frac12\int J\cdot\nabla|w|^2=0.
 \label{eq:cancellation}
\end{align}
Testing the full-space equation by $w$ and taking real parts, we have 
$\int\wpi K\nabla w\cdot\nabla\overline w=0$. Hence $w$ is constant,
proving the kernel statement and normalized uniqueness.
\end{proof}

For $h\in H$, set $w=h/\wpi$,
$\mathcal E=\int_\Omega\wpi|\nabla w|^2\dd V$, 
and $\mathcal N=\int_\Omega\wpi|w|^2\dd V$.
The quotient turns the stationary mode into a constant, while $wJ$ has
zero real energy by \eqref{eq:cancellation}.
Since $\int\wpi w=0$, weighted variance and Poincar\'e's inequality give
\begin{equation}
 \mathcal N\le C_\wpi^2\mathcal E,
 \quad C_\wpi^2=
 \frac{\wpi_{\max}}{\lambda_1^{\mathrm{Neu}}\wpi_{\min}},
 \label{eq:weightedP}
\end{equation}
where $\lambda_1^{\mathrm{Neu}}>0$ is the first positive Neumann eigenvalue. 
Indeed, since $\overline w_\Omega$ is the ordinary mean,
$\mathcal N\le\int\wpi|w-\overline w_\Omega|^2
\le\wpi_{\max}\|\nabla w\|_2^2/\lambda_1^{\mathrm{Neu}}\le C_\wpi^2\mathcal E$.
In particular, $\mathcal E>0$ for $h\ne0$ in $H$.

Applying \eqref{eq:weakA} to an eigenfunction $h$ satisfying
$A_\sigma h=\lambda h$ gives
\begin{equation}
 -\diver j_h=\lambda(\sigma I-\Delta)h,\quad
 \nu\cdot j_h=\lambda\partial_\nu h,\quad\int_\Omega h=0.
 \label{eq:pencil}
\end{equation}
Set $G=\norm{\nabla\log\wpi}_{L^\infty}$ and
$V_J=\norm{J/\wpi}_{L^\infty}$. 
Here $G$ measures variation of the weight, and $V_J$ measures its steady current. 
The following choices are sufficient for Theorem~\ref{thm:continuous}:
\begin{equation}
 \sigma_*=\frac{G^2}{2}+\frac{V_JG}{\kappa_-},\quad
 M^2=\kappa_+^2+V_J^2C_\wpi^2,\quad
 \tau_*=\frac{\kappa_-}{M^2}.
 \label{eq:constants2}
\end{equation}

\begin{lemma}[Eigenvalue estimate]\label{lem:quotient}
Assume \eqref{eq:parameters} and set
$D_\sigma=1+GC_\wpi+\sigma C_\wpi^2$. Every eigenvalue $\lambda$ of
$A_\sigma$ on either $H$ or $X$ satisfies
\begin{equation}
 |1-\tau\lambda|\le\qe:=
 \left(1-\frac{\tau(\kappa_--\tau M^2)}{D_\sigma^2}\right)^{1/2}<1.
 \label{eq:qe}
\end{equation}
\end{lemma}
\begin{proof}
For an eigenfunction $h\in H$, set $w=h/\wpi$ and test
\eqref{eq:weakA} by $\phi=w-\mean w\in H$.
Since $\nabla\phi=\nabla w$ and $\int_\Omega h\dd V=0$,
subtracting the ordinary mean leaves both sides of the weak equation unchanged.
Substitution of $h=\wpi w$ gives $\mathfrak a=\lambda\mathfrak b$, where
\begin{align}
 \mathfrak a&=\int\wpi K\nabla w\cdot\nabla\overline w
                      +\int wJ\cdot\nabla\overline w,\label{eq:ab1}\\
 \mathfrak b&=\mathcal E+\sigma\mathcal N
                      +\int w\nabla\wpi\cdot\nabla\overline w.
 \label{eq:ab2}
\end{align}
The cancellation \eqref{eq:cancellation} and Cauchy--Schwarz give
\begin{align*}
 \kappa_-\mathcal E\le\Rea\mathfrak a\le\kappa_+\mathcal E,
 &\quad |\Ima\mathfrak a|\le V_J\sqrt{\mathcal N\mathcal E},\\
 \Rea\mathfrak b\ge\mathcal E+\sigma\mathcal N
                         -G\sqrt{\mathcal N\mathcal E},
 &\quad |\Ima\mathfrak b|\le G\sqrt{\mathcal N\mathcal E}.
\end{align*}
Young's inequality implies
$\Rea\mathfrak b\ge\tfrac12\mathcal E+
(\sigma-G^2/2)\mathcal N>0$. 
Since $\frac{\mathfrak a}{\mathfrak b} = 
(\mathfrak a \overline{\mathfrak b}) / |\mathfrak b|^2$, 
the sign of $\Rea \lambda$ is determined by $\Rea(\mathfrak a\overline{\mathfrak b})$. 
The real part satisfies
\begin{align}
 \Rea(\mathfrak a\overline{\mathfrak b})
 \ge\tfrac12(\Rea\mathfrak a)\mathcal E+
       \left[\kappa_-(\sigma-G^2/2)-V_JG\right]\mathcal N\mathcal E 
 \ge\tfrac{\kappa_-}{2}\mathcal E^2.
 \label{eq:realquotient}
\end{align}
The last inequality uses $\sigma\ge\sigma_*$. Estimate~\eqref{eq:weightedP} also gives
$ |\mathfrak a|^2\le M^2\mathcal E^2 $ and 
$ |\mathfrak b|\le D_\sigma\mathcal E$, hence
\[
 |1-\tau\lambda|^2
 =1-\frac{2\tau\Rea(\mathfrak a\overline{\mathfrak b})
                     -\tau^2|\mathfrak a|^2}{|\mathfrak b|^2}
 \le1-\frac{\tau(\kappa_--\tau M^2)}{D_\sigma^2}.
\]
This weak argument covers complex eigenvalues as well.
\end{proof}

\subsection{Full spectral analysis}

An infinite-dimensional operator may have spectral points that are not eigenvalues, 
so Lemma~\ref{lem:quotient} alone does not control the full spectrum.
We locate the remaining part through diffusion and treat the drift as a compact
perturbation.

\paragraph{Energy space.}
Split $A_\sigma=B_\sigma+\Dcal_\sigma$ on $H$ into diffusion and drift parts by
\begin{equation}
 p_\sigma(B_\sigma h,\phi)=\int K\nabla h\cdot\nabla\overline\phi, \quad
 p_\sigma(\Dcal_\sigma h,\phi)=-\int hKd\cdot\nabla\overline\phi.
 \label{eq:Bdef}
\end{equation} 
Lax-Milgram guarantees that both operators are bounded 
with respect to the norm $\norm{\cdot}_\sigma$ defined in \eqref{eq:energy}.

\begin{lemma}[Energy decomposition]\label{lem:energydecomp}
The operator $B_\sigma$ is positive and self-adjoint in $p_\sigma$, with
\begin{equation}
 \kappa_\sigma I\preceq B_\sigma\preceq\kappa_+I,
 \quad
 \kappa_\sigma=
 \frac{\kappa_-\lambda_1^{\mathrm{Neu}}}{\lambda_1^{\mathrm{Neu}}+\sigma}>0.
 \label{eq:Bbounds}
\end{equation}
The operator $\Dcal_\sigma$ is compact on $H$. Consequently
\begin{equation}
 \Tcal=\Tcal_0+\Kcal,\quad
 \Tcal_0=I-\tau B_\sigma,\quad \Kcal=-\tau\Dcal_\sigma,
 \label{eq:split}
\end{equation}
where, for $0<\tau<2/\kappa_+$,
\begin{equation}
 \norm{\Tcal_0}_{\sigma\to\sigma}\le
  q_B:=\max\{|1-\tau\kappa_\sigma|,|1-\tau\kappa_+|\}<1.
 \label{eq:qB}
\end{equation}
\end{lemma}
\begin{proof}
Symmetry and the upper bound follow from \eqref{eq:Bdef}. Poincar\'e's
inequality gives
$p_\sigma(h,h)\le(1+\sigma/\lambda_1^{\mathrm{Neu}})\norm{\nabla h}_2^2$,
which proves the lower bound. Definition of the drift term gives
\[
 |p_\sigma(\Dcal_\sigma h,\phi)|
 \le\|Kd\|_\infty\|h\|_2\norm\phi_\sigma.
\]
Thus $\Dcal_\sigma$ factors through the compact embedding
$H^1\hookrightarrow L^2$ and the bounded Riesz inverse for $p_\sigma$.
The self-adjoint spectral theorem proves \eqref{eq:qB}.
\end{proof}

For $z\notin[\kappa_\sigma,\kappa_+]$, factor
\[
 A_\sigma-zI=
 \bigl[I+\Dcal_\sigma(B_\sigma-zI)^{-1}\bigr](B_\sigma-zI).
\]
The bracket is identity plus a compact analytic family, and invertible for
large $|z|$. Analytic Fredholm theory \cite{Kato} implies there are only isolated
eigenvalues for $A_\sigma$ outside the interval. Combining \eqref{eq:qB},
Lemma~\ref{lem:quotient}, and spectral mapping gives
\begin{equation}
 \rad_H(\Tcal)\le\max\{ q_B,\qe\}<1.
 \label{eq:Hrad}
\end{equation}

\paragraph{Strong space.}
Since the nonlinear map acts on $X$, we must also control the spectrum
there. Introduce the comparison operator $A_{\rm pr}$ by
\begin{equation}
 \int\nabla A_{\rm pr}h\cdot\nabla\overline\phi
       =\int K\nabla h\cdot\nabla\overline\phi,
 \quad \phi\in H.
 \label{eq:Apr}
\end{equation}
It is self-adjoint in the gradient inner product on $H$ and bounded on $X$ by
Lemma~\ref{lem:solver}.

\begin{lemma}[Spectrum in strong space]\label{lem:conormal}
The spectrum of $A_{\rm pr}$ on $X$ is contained in
$[\kappa_-,\kappa_+]$.
\end{lemma}
\begin{proof}
For $z\notin I_K=[\kappa_-,\kappa_+]$ and $w\in X$, the energy-space
spectral theorem gives a unique $h\in H$ with $(A_{\rm pr}-zI)h=w$ and 
$\norm h_{H^1}\le C_z\norm w_{H^1}$. Since constants contribute zero, 
the weak equation holds for all $H^1$ tests and reads
\begin{equation}
 -\diver((K-zI)\nabla h)=-\Delta w,
 \quad \nu\cdot(K-zI)\nabla h=\partial_\nu w.
 \label{eq:resolventpde}
\end{equation}

Set $B=K-zI$. Since the segment $I_K-z$ avoids zero, choosing $\theta$
opposite to the argument of its nearest point to zero gives
\[
 \Rea(e^{i\theta}B\zeta\cdot\overline\zeta)
 \ge\delta_z|\zeta|^2,
 \quad \delta_z=\dist(z,I_K)>0.
\]
Following the half-space test in \cite[Section~2.2]{CostabelDaugeNicaise},
choose boundary normal coordinates $(s,t)$ in a neighborhood of fixed small
normal width $\varepsilon>0$, with $t=0$ on the boundary and
$0<t<\varepsilon$ inside. Freeze the transformed principal matrix as
$\widetilde B$. The real coordinate change preserves rotated coercivity.
The frozen homogeneous conormal model is posed on
$\mathbb R^{n-1}\times(0,\infty)$, where $n=\dim\Omega$.
For a decaying mode $e^{i\xi\cdot s}q(t)$, $\xi\ne0$, integrate by parts
on $(0,T)$, where $T>0$ is an arbitrary truncation of the model half-line. 
Letting $T\to\infty$ gives
\[
 0=\Rea\int_0^\infty e^{i\theta}\widetilde B(i\xi q,q')
                    \cdot\overline{(i\xi q,q')}\dd t
 \ge c_z\int_0^\infty\bigl(|\xi|^2|q|^2+|q'|^2\bigr)\dd t, 
 \quad c_z > 0.
\]
The endpoint terms vanish by zero conormal data, 
thus $q=0$. Along the coercive path
$(1-\lambda)e^{i\theta}\widetilde B+\lambda I$, $0\le\lambda\le1$,
the normal leading coefficient stays nonzero and no root crosses the real
axis. At $I$, the normal roots are $\pm i|\xi|$. Hence exactly one root
lies in the upper half-plane throughout the path.
It verifies proper ellipticity and the complementing condition.

Since $w\in X$ and $K\in C^{1,1}\subset C^{1,\alpha}$, 
boundary Sobolev and Schauder regularity \cite{ADN1,ADN2} imply
$h\in C^{2,\alpha}(\overline\Omega)$. Their estimate and the $H^1$ bound yield
\[
 \norm h_X\le C_z\bigl(
 \norm{\Delta w}_{C^{0,\alpha}}+
 \norm{\partial_\nu w}_{C^{1,\alpha}(\partial\Omega)}+
 \norm h_{H^1}\bigr)
 \le C'_z\norm w_X.
\]
Zero mean and uniqueness follow from the weak problem, 
so $(A_{\rm pr}-zI)^{-1}$ is bounded on $X$. 
\end{proof}

For $h\in X$, the drift term $\mathcal D_\sigma h$
is equivalently the zero-mean solution $w_h$ of
\[
 (\sigma I-\Delta)w_h=\diver(Kdh),\quad
 \partial_\nu w_h=-\nu\cdot Kdh.
\]
The compatibility constant is zero by the divergence theorem. Let 
$R_\sigma^{\mathrm{Neu}}$ be the zero-mean homogeneous Neumann inverse.
Then $A_\sigma$ can be decomposed on $X$ as
\begin{equation}
 A_\sigma=A_{\rm pr}-\sigma R_\sigma^{\mathrm{Neu}}A_{\rm pr}
                       +\Dcal_\sigma.
 \label{eq:strongsplit}
\end{equation}
The bounded map $A_{\rm pr}:X\to X$ and Neumann regularity give
$R_\sigma^{\mathrm{Neu}}A_{\rm pr}h\in C^{4,\alpha}$, so the shift term
is compact on $X$. For the drift, the inclusion
$X\hookrightarrow C^{1,\alpha}(\overline\Omega)$ is compact, and $Kd\in C^{1,\alpha}$
makes $h\mapsto(\diver(Kdh),-\nu\cdot Kdh)$ bounded into
$C^{0,\alpha}(\overline\Omega)\times C^{1,\alpha}(\partial\Omega)$.
These data satisfy Neumann compatibility by the divergence theorem.
Lemma~\ref{lem:solver} maps them boundedly into $X$, hence this factorization
makes $\Dcal_\sigma$ compact on $X$. Thus $A_\sigma-A_{\rm pr}$ is compact.

Another Fredholm argument and Lemma~\ref{lem:quotient} show that every
spectral point of $A_\sigma$ on $X$ outside $[\kappa_-,\kappa_+]$ is an eigenvalue.
Under \eqref{eq:parameters} and \eqref{eq:constants2}, we have $\tau\kappa_+<1$.
The interval therefore contributes at most $1-\tau\kappa_-$ to the spectral
radius of $\Tcal$, giving
\begin{equation}
 \rad_X(\Tcal)\le\max\{1-\tau\kappa_-,\qe\}<1.
 \label{eq:Xrad}
\end{equation}

If $d=0$, $A_\sigma=B_\sigma$ has energy spectrum in
$[\kappa_\sigma,\kappa_+]$. Spectral points of $X$ outside
$[\kappa_-,\kappa_+]$ are eigenvalues in $H$, and the comparison interval
lies in $[\kappa_\sigma,\kappa_+]$. Hence under \eqref{eq:driftfree}, 
\begin{equation}
 \rad_H(\Tcal),\ \rad_X(\Tcal)\le q_B<1.
 \label{eq:driftfreespectrum}
\end{equation}
Thus zero reference drift permits the unshifted mean-free Poisson inverse.

\subsection{Completion of the proof}
\label{sec:continuousnonlinear}
\begin{proof}[Proof of Theorem~\ref{thm:continuous}]
Since $u^\ast$ solves \eqref{eq:physicalproblem},
$F(u^\ast)=\Psi(u^\ast)=m(u^\ast)=0$. Uniqueness of the normalized
Neumann solution gives $\varphi(u^\ast)=0$, hence $S(u^\ast)=u^\ast$.

Equations~\eqref{eq:Xrad} and~\eqref{eq:driftfreespectrum} give
$\rad_X(\Tcal)<1$. Choose
$\rad_X(\Tcal)<s<\gamma<q<1$. The spectral-radius formula implies
$\norm{\Tcal^n}_{X\to X}\le Cs^n$ for some $C\ge1$. Hence
\begin{equation}
 \norm h_\dagger=\sup_{n\ge0}\gamma^{-n}\norm{\Tcal^nh}_X
 \label{eq:adapted}
\end{equation}
defines an equivalent norm with
$\norm h_X\le\norm h_\dagger\le C\norm h_X$
and $\norm{\Tcal h}_\dagger\le\gamma\norm h_\dagger$. 
By Lemma~\ref{lem:solver} and norm equivalence, there is a
small closed ball $\mathcal B\subset\Ucal$ about $u^\ast$ with
$\norm{DS(u)-\Tcal}_{\dagger\to\dagger}<q-\gamma$ for all $u\in \mathcal{B}$, 
hence $\sup_{u\in\mathcal{B}}\norm{DS(u)}_{\dagger\to\dagger}\le q$.
Integration along segments makes $S$ a $q$-contraction on this ball.
Since $S(u^\ast)=u^\ast$,
\[
 \norm{S(u)-u^\ast}_\dagger\le q\norm{u-u^\ast}_\dagger,
 \quad u\in\mathcal B.
\]
Thus $\mathcal B$ is invariant and complete in the equivalent norm.
Banach's theorem identifies $u^\ast$ as the unique fixed point and gives
$\norm{u^n-u^\ast}_\dagger\le q^n\norm{u^0-u^\ast}_\dagger$ by iteration.
Choose $\delta$ so that the $X$-ball of radius $\delta$ lies in
$\mathcal B$, and norm equivalence proves \eqref{eq:continuousconvergence}. 
The local boundedness of the derivatives of $F,\Psi,m$ yields \eqref{eq:continuousresiduals}.
The update preserves the zero mean as well.
\end{proof}

\label{page:continuous-end}
\section{Spectral Galerkin approximation}
\label{sec:spectral}
We discretize the continuous correction \eqref{eq:continuousmap} on
spherical caps using Fourier--Jacobi modes on a reference disk.

\subsection{The discrete Poisson correction and algorithm}
For $N\ge1$, let $V_N\subset X$ be the zero-mean spectral space of order $N$.
Given $u_N^n\in V_N$ and fixed $\sigma,\tau,\Psi$, find $v_N^n\in V_N$ satisfying
\begin{equation}
 p_\sigma(v_N^n,\phi_N)
 =\int_\Omega F(u_N^n)\phi_N\dd V
  +\int_{\partial\Omega}\Psi(u_N^n)\phi_N\dd S,
 \quad \phi_N\in V_N.
 \label{eq:galerkin}
\end{equation}
The zero-mean test functions eliminate the compatibility term.

Let $d_N = \dim V_N$. 
In a real basis $\{\phi_j\}_{j=1}^{d_N}$, write
$u_N=\sum_jU_j\phi_j$ and define
\begin{align}
 (L_N)_{ij}&=\int_\Omega\nabla\phi_j\cdot\nabla\phi_i\dd V, \quad
 (M_N)_{ij}=\int_\Omega\phi_j\phi_i\dd V,\label{eq:PN}\\
 (R_N(U))_i&=\int_\Omega F(u_N)\phi_i\dd V
 +\int_{\partial\Omega}\Psi(u_N)\phi_i\dd S.\nonumber
\end{align}
One step can be expressed as 
\begin{equation}
 (L_N+\sigma M_N)V^n=R_N(U^n),
 \quad U^{n+1}=U^n-\tau V^n.
 \label{eq:matrixiteration}
\end{equation}
The fixed matrix is positive definite even at $\sigma=0$ and is factored once.
The boundary load supplies the Neumann data~\cite{Atkinson}.
Algorithm~\ref{alg:fixedmatrix} combines quadrature norms of $F$ and
$a^{-1}\Psi=\psi\circ T_u$, with a fixed boundary scale $\ell_\partial>0$.
The boundary term is omitted for exactly imposed boundary conditions.

\begin{algorithm}[H]
\caption{Spectral Galerkin Poisson iteration}\label{alg:fixedmatrix}
\begin{algorithmic}[1]
\Require $N,\sigma,\tau$, frozen weight $a$, boundary scale $\ell_\partial>0$, $U^0$, budget $n_{\max}$,
$\epsilon_{\rm abs}>0$, $n_{\rm check}\ge1$, $\epsilon_{\rm dec}>0$.
\State Construct the basis and factor $P_N=L_N+\sigma M_N$.
\For{$n=0,\ldots,n_{\max}$}
 \State Evaluate $u_N^n$, its derivatives, $T_{u_N^n}$, $F(u_N^n)$ and $\Psi(u_N^n)$ where defined.
 \State $e_{F,n}\gets\|F(u_N^n)\|_{L^2(\Omega)}/\|f\|_{L^2(\Omega)}$.
 \State $e_{\partial,n}\gets\|a^{-1}\Psi(u_N^n)\|_{L^2(\partial\Omega)}/|\partial\Omega|^{1/2}$.
 \State $\eta_n\gets\sqrt{e_{F,n}^2+(e_{\partial,n}/\ell_\partial)^2}$.
 \State \textbf{if} evaluation fails or produces nonfinite values, \textbf{return} with the failure reason.
 \State \textbf{if} $\eta_n\le\epsilon_{\rm abs}$, \textbf{return} $u_N^n$ (absolute tolerance).
 \If{$n\ge n_{\rm check}$}
  \State $d_n\gets(\eta_{n-1}-\eta_n)/\eta_{n-1}$.
  \State \textbf{if} $d_n<0$, \textbf{return} $u_N^n$ (residual increase).
  \State \textbf{if} $d_n<\epsilon_{\rm dec}$, \textbf{return} $u_N^n$ (insufficient decrease).
 \EndIf
 \State \textbf{if} $n=n_{\max}$, \textbf{return} $u_N^n$ (iteration budget).
 \State Assemble $R_N(U^n)$ from signed residuals and solve $P_NV^n=R_N(U^n)$.
 \State Update $U^{n+1}=U^n-\tau V^n$.
\EndFor
\end{algorithmic}
\end{algorithm}

The energy projection $\Pi_N:H\to V_N$ satisfies
$p_\sigma(\Pi_Nh-h,\phi_N)=0$ for all $\phi_N\in V_N$.
Thus $v_N=\Pi_N\varphi(u_N)$ and
\begin{equation}
 S_N(u_N)=u_N-\tau\Pi_N\varphi(u_N),
 \quad u_N\in V_N\cap\Ucal.
 \label{eq:SN}
\end{equation}

\subsection{Fourier--Jacobi discretization on a spherical cap}
\label{sec:disk_spectral}
Let $\Omega=\{x\in\mathbb S^2:x\cdot e_3>\cos\theta_0\}$,
$0<\theta_0<\pi$, and $D=\{\xi\in\bR^2:|\xi|<1\}$.
Set $r_c=\tan(\theta_0/2)$ and $\varrho=|\xi|$.
The stereographic parametrization
$\Phi:\overline D\to\overline\Omega$ and its conformal factor are
\[
 \Phi(\xi)=\frac{(2r_c\xi_1,\,2r_c\xi_2,\,1-r_c^2\varrho^2)}
                    {1+r_c^2\varrho^2},\quad
 \omega(\varrho)=\frac{2r_c}{1+r_c^2\varrho^2}.
\]
Write $\widehat h=h\circ\Phi$. The conformal metric is $\omega^2I$,
so $\dd V=\omega^2\dd A$ and the metric factors cancel in the
two-dimensional gradient energy. Using Euclidean derivatives on $D$,
\begin{equation}
 p_\sigma(w,\phi)=
 \int_D\nabla\widehat w\cdot\nabla\overline{\widehat\phi}\dd A
 +\sigma\int_D\omega^2\widehat w\,\overline{\widehat\phi}\dd A,
 \quad\dd S=\omega(1)\dd\theta.
 \label{eq:diskenergy}
\end{equation}
Let $\Pcal_N(D)$ denote polynomials of total degree at most $N$ and set
\begin{equation}
 V_N=\left\{h:\widehat h\in\Pcal_N(D),\quad
               \int_D\omega^2\widehat h\dd A=0\right\},
 \quad d_N=\frac{(N+1)(N+2)}2-1.
 \label{eq:VN}
\end{equation}
In disk polar coordinates $(\varrho,\theta)$, real Fourier--Jacobi modes are
\begin{equation}
 \varrho^mP_k^{(0,m)}(2\varrho^2-1)\cos(m\theta),\quad
 \varrho^mP_k^{(0,m)}(2\varrho^2-1)\sin(m\theta),
 \quad 2k+m\le N.
 \label{eq:basis}
\end{equation}
The constant and the zero sine modes are excluded here. 
To enforce zero mean on the spherical cap, subtract the weighted mean
from each mode $Z_j$:
\begin{equation}
 \widehat\phi_j=Z_j-
    \frac{\int_D\omega^2 Z_j\dd A}{\int_D\omega^2\dd A}.
 \label{eq:meanbasis}
\end{equation}
In the implementation, coefficients are computed using
periodic trapezoidal quadrature and FFT in $\theta$, and
Gauss--Legendre quadrature in $s=\varrho^2$ with precomputed radial Jacobi
values and derivatives.

\begin{remark}[Square domains]\label{rem:square_spectral}
The same spectral discretization scheme can be implemented on a square using a
cosine basis, as in the refractor experiment in
Section~\ref{sec:square_refractor}. The convergence analysis in this
paper is restricted to spherical caps.
\end{remark}

\section{Local convergence of the Galerkin iteration}
\label{sec:discrete}
In this section, we extend the local convergence theorem to the Galerkin discretization. $\Omega$ is the spherical cap of
Section~\ref{sec:disk_spectral}, with its fixed disk parametrization.
All discrete statements refer to these Fourier--Jacobi spaces.
The analysis concerns the exact Galerkin map \eqref{eq:SN}.
\subsection{The discrete local convergence theorem}
For the strong norm on zero-mean differences, choose
\begin{equation}
 u_N^\ast\in\mathop{\rm argmin}_{w_N\in V_N}
                     \norm{w_N-u^\ast}_X,
 \quad
 \epsN=\norm{u_N^\ast-u^\ast}_X.
 \label{eq:bestapprox}
\end{equation}
Smooth functions are dense in $H^1$. Disk polynomials can
approximate smooth functions in $H^4$ \cite{Figueroa}.
By the compact embedding theorem, 
the spectral spaces have two basic approximation properties:
\begin{equation}
 \Pi_N \longrightarrow I \text{ in }H,
 \quad \epsN\longrightarrow0.
 \label{eq:qualitativeapprox}
\end{equation}
Since $V_N$ is finite-dimensional, its weak and strong norms are equivalent.
Define the finite inverse constant by
\begin{equation}
 \invN=\max\left\{1,\sup_{0\ne h_N\in V_N}
                    \frac{\norm{h_N}_X}{\norm{h_N}_\sigma}\right\}.
 \label{eq:inverseconstant}
\end{equation}
\begin{theorem}[Discrete local convergence]\label{thm:discrete}
Assume the hypotheses and either parameter range of
Theorem~\ref{thm:continuous}. Use the spaces \eqref{eq:VN} for a fixed smooth
disk parametrization, the same fixed $\sigma,\tau$ and matched boundary
weight. There exist
$N_0$, $\rho>0$, $C,C_a\ge1$, $q\in(0,1)$, and norms
\[
 \norm{h_N}_\sigma\le\norm{h_N}_{\dagger,N}
 \le C_a\norm{h_N}_\sigma.
\]
When $N\ge N_0$, $S_N$ is a
$q$-contraction on the closed invariant ball
\begin{equation}
 \mathcal B_N=
 \{u_N\in V_N:\norm{u_N-u_N^\ast}_{\dagger,N}\le r_N\},
 \quad r_N=\rho/\invN^2.
 \label{eq:discreteball}
\end{equation}
This ball lies in $\Ucal$, where $\mathcal U$ is the regular neighborhood from
Lemma~\ref{lem:neighborhood}. It contains a unique fixed point $u_N^\infty$,
with
\begin{equation}
 \norm{u_N^\infty-u_N^\ast}_{\dagger,N}\le C\epsN,
 \quad
 \norm{u_N^\infty-u^\ast}_X\le C(1+\invN)\epsN.
 \label{eq:discreteerrorX}
\end{equation}
Every $u_N^0\in\mathcal B_N$ generates iterates staying in that ball, and
\begin{align}
 \norm{u_N^n-u^\ast}_X
 &\le C\invN q^n\norm{u_N^0-u_N^\infty}_\sigma
       +C(1+\invN)\epsN.
 \label{eq:totalerror}
\end{align}
The strong residuals of the fixed point satisfy
\begin{equation}
 \norm{F(u_N^\infty)}_{C^{0,\alpha}}
 +\norm{\Psi(u_N^\infty)}_{C^{1,\alpha}(\partial\Omega)}
 +|m(u_N^\infty)|
 \le C(1+\invN)\epsN.
 \label{eq:discreteresidual}
\end{equation}
\end{theorem}

\subsection{Uniform stability of the Galerkin linearizations}
Differentiate the residual at $u^*$. 
The boundary--flux identity gives the compression of the continuous
linearized update $\Tcal=DS(u^\ast)$ as
\begin{align}
 (C_N)_{ij}
 &=p_\sigma(A_\sigma\phi_j,\phi_i)
 =\int_\Omega K(\nabla\phi_j-d\phi_j)\cdot\nabla\phi_i\dd V,
 \nonumber\\
 \widehat\Tcal_N
 &=\Pi_N\Tcal|_{V_N}
 =I-\tau(L_N+\sigma M_N)^{-1}C_N.
 \label{eq:CN}
\end{align}

\begin{lemma}[Uniform spectral bounds]\label{lem:galerkin}
Under either parameter range of Theorem~\ref{thm:continuous}, for every
\begin{equation}
 \max\{ q_B,\rad_H(\Tcal)\}<q_0<1
 \label{eq:bridgeq}
\end{equation}
there exist $N_0,C_0$ such that, for all $N\ge N_0$ and $k\ge0$,
\begin{equation}
 \spec(\widehat\Tcal_N)\subset\{|z|<q_0\},
 \quad
 \norm{\widehat\Tcal_N^{\,k}}_{\sigma\to\sigma}
 \le C_0q_0^k.
 \label{eq:uniformpowers}
\end{equation}
If $d=0$, every nontrivial trial space already satisfies
$\norm{\widehat\Tcal_N}_{\sigma\to\sigma}\le q_B<1$.
\end{lemma}
\begin{proof}
All operator norms in this proof are energy norms $\norm{\cdot}_{\sigma\to\sigma}$.
The continuous analysis gives $\Tcal=\Tcal_0+\Kcal$ on the complexified
energy space, with $\norm{\Tcal_0}\le q_B<1$, $\Kcal$ compact, and
$\rad_H(\Tcal)<1$. For $\Tcal_{0,N}=\Pi_N\Tcal_0|_{V_N}$ and $|z| > q_B$,
orthogonality and the Neumann series imply
\begin{equation}
 \norm{(zI-\Tcal_{0,N})^{-1}}
 \le(|z|- q_B)^{-1}.
 \label{eq:baseinverse}
\end{equation}
By \eqref{eq:qualitativeapprox}, each fixed power converges strongly, and 
the geometric tails are uniformly bounded. Hence, for $y\in H$,
\begin{equation}
 (zI-\Tcal_{0,N})^{-1}\Pi_Ny
 \longrightarrow(zI-\Tcal_0)^{-1}y
 \label{eq:strongres}
\end{equation}
locally uniformly for $|z|> q_B$.
If $zI-\widehat\Tcal_N$ is not uniformly invertible on $q_0\le|z|\le R$, where
$R>\norm{\Tcal}+1$, there would be $N_j\to\infty$, $z_j\to z$, and 
$h_j\in V_{N_j}$ with $\norm{h_j}_\sigma = 1$, such that
\begin{equation}
 (z_jI-\widehat\Tcal_{N_j})h_j\longrightarrow0 .
 \label{eq:approxkernel}
\end{equation}
Compactness gives a subsequence $\Kcal h_j\to y$, and
\eqref{eq:baseinverse}, \eqref{eq:strongres} imply
\[
 h_j=(z_jI-\Tcal_{0,N_j})^{-1}\Pi_{N_j}\Kcal h_j+o(1)
 \longrightarrow h=(zI-\Tcal_0)^{-1}y.
\]
Thus $\norm h_\sigma=1$ and $\Tcal h=zh$, contradicting
$|z|\ge q_0>\rad_H(\Tcal)$. For $|z|>R$, the bound
$\norm{\widehat\Tcal_N}\le\norm{\Tcal}$ gives invertibility directly.
Let $M$ be the uniform resolvent bound on $|z| = q_0$. By Cauchy's formula,
\[
 \widehat\Tcal_N^{\,k}=\frac1{2\pi i}\int_{|z|=q_0}
 z^k(zI-\widehat\Tcal_N)^{-1}\dd z.
\]
It implies $\norm{\widehat\Tcal_N^{\,k}}\le Mq_0^{k+1}$. 
Taking $C_0=\max\{1,Mq_0\}$ proves \eqref{eq:uniformpowers} uniformly in $N$ and $k$.
If $d=0$, $\Kcal=0$. Together with orthogonal
compression, this gives the stated one-step bound.
\end{proof}

Choose $q_0<\gamma<q<1$ and define
\begin{equation}
 \norm{h_N}_{\dagger,N}=\sup_{k\ge0}\gamma^{-k}
 \norm{\widehat\Tcal_N^{\,k}h_N}_\sigma.
 \label{eq:discretenorm}
\end{equation}
Lemma~\ref{lem:galerkin} gives a constant $C_a$ independent of large $N$ such that 
\begin{equation}
 \norm{h_N}_\sigma\le\norm{h_N}_{\dagger,N}\le C_a\norm{h_N}_\sigma,
 \quad\norm{\widehat\Tcal_N}_{\dagger,N}\le\gamma.
 \label{eq:uniformnorm}
\end{equation}

\subsection{Nonlinear estimates}
We estimate the derivative perturbation using the approximation error
and the distance from the best approximant $u_N^\ast$.

\begin{lemma}[Nonlinear derivative bounds]
\label{lem:discretebounds}
Shrink $\mathcal U$ from Lemma~\ref{lem:neighborhood}. 
There exist constants $C_Y,C_Z$, independent of sufficiently large $N$,
such that
\begin{equation}
 \norm{S_N(u_N^\ast)-u_N^\ast}_{\dagger,N}\le C_Y\epsN.
 \label{eq:defect}
\end{equation}
On any closed ball centered at $u_N^\ast$, of radius $r$ in the adapted
norm and contained in $\mathcal U$,
\begin{equation}
 \sup_{\norm{u_N-u_N^\ast}_{\dagger,N}\le r}
 \norm{DS_N(u_N)-\widehat\Tcal_N}_{\dagger,N}
 \le C_Z\invN(\epsN+\invN r).
 \label{eq:derivativeconsistency}
\end{equation}
\end{lemma}
\begin{proof}
Projection stability gives $\norm{\Pi_Nz}_{\dagger,N}\le C_a\norm z_\sigma$.
By Lemma~\ref{lem:solver}, the derivative map
$u\mapsto D\varphi(u)\in\mathcal L(X,H)$ has a uniform bound $M_H$
and a Lipschitz constant $L_H$ on $\mathcal U$. Since $\varphi(u^\ast)=0$,
\[
 \norm{S_N(u_N^\ast)-u_N^\ast}_{\dagger,N} = \tau \norm{\Pi_N\varphi(u_N^\ast)}_{\dagger,N}
 \le C_a\tau M_H\norm{u_N^\ast-u^\ast}_X=C_Y\epsN.
\]
For $h_N\in V_N$,
$(DS_N(u_N)-\widehat\Tcal_N)h_N
=-\tau\Pi_N(D\varphi(u_N)-D\varphi(u^\ast))h_N$. Hence
\[
 \norm{(DS_N(u_N)-\widehat\Tcal_N)h_N}_{\dagger,N}
 \le C_a\tau L_H\norm{u_N-u^\ast}_X\invN\norm{h_N}_{\dagger,N}.
\]
Set $C_Z=C_a\tau L_H$ and use 
$\norm{u_N-u^\ast}_X\le \norm{u_N - u_N^*}_X + \norm{u_N^* - u_*}_X
\le\epsN+\invN r$. 
Then we prove \eqref{eq:derivativeconsistency}.
\end{proof}

\Needspace{8\baselineskip}
\subsection{Proof of the discrete convergence theorem}
\begin{lemma}[Polynomial approximation and inverse bounds]
\label{lem:polynomials}
For the spaces \eqref{eq:VN} and fixed $0<\alpha<1$, the inverse constant satisfies
\begin{equation}
 1\le\invN\le C_\alpha(N+1)^6.
 \label{eq:inversepolynomial}
\end{equation}
For $\widehat u^\ast=u^\ast\circ\Phi\in H^\ell(D)$, with integer
$\ell\ge4$,
\begin{equation}
 \epsN\le C_\ell(N+1)^{15/2-\ell}
                    \norm{\widehat u^\ast}_{H^\ell(D)}.
 \label{eq:bestapproxrate}
\end{equation}
In particular, the smooth reference of Theorem~\ref{thm:continuous} satisfies
\begin{equation}
 \invN^2\epsN\longrightarrow0,
 \quad\invN\epsN\longrightarrow0.
 \label{eq:approxinverse}
\end{equation}
\end{lemma}
\begin{proof}
In the fixed chart, Poincar\'e gives
$\norm{\widehat h}_{H^1(D)}\le C\norm h_\sigma$, also for $\sigma=0$.
Applying the Markov estimate
$\norm{\nabla p}_{L^2}\le C(N+1)^2\norm p_{L^2}$
\cite[Lemma~2.11]{Figueroa} three times gives
$\norm p_{H^4}\le C(N+1)^6\norm p_{H^1}$.
Then $H^4\hookrightarrow C^{2,\alpha}$ implies \eqref{eq:inversepolynomial}.

For the projection $Q_N: L^2(D)\to \mathcal{P}_N(D)$, the unweighted
estimate \cite[Theorem~3.9]{Figueroa} at derivative order four gives
\begin{equation}
 \norm{w-Q_Nw}_{H^4(D)}
 \le C_\ell(N+1)^{15/2-\ell}\norm w_{H^\ell(D)}.
 \label{eq:unitprojection}
\end{equation}
Construct the explicit mean-corrected approximant
\begin{equation}
 \widehat{\widetilde u}_N=Q_N\widehat u^\ast+c_N,
 \quad
 c_N=\frac{\int_D\omega^2(\widehat u^\ast-Q_N\widehat u^\ast)\dd A}
              {\int_D\omega^2\dd A}.
 \label{eq:meanapproximation}
\end{equation}
Then $\widetilde u_N\in V_N$ and
$|c_N|\le C\norm{\widehat u^\ast-Q_N\widehat u^\ast}_{L^2}$.
Embedding $H^4\hookrightarrow C^{2,\alpha}$ and 
the minimum definition of $\epsN$ imply 
$\epsN\le C\norm{\widehat u^\ast-Q_N\widehat u^\ast}_{H^4}$,
which proves \eqref{eq:bestapproxrate}. Consequently,
\begin{equation}
 \invN^2\epsN
 \le C_\ell(N+1)^{39/2-\ell}
                \norm{\widehat u^\ast}_{H^\ell(D)}\longrightarrow0
 \quad \text{when } \ell>39/2.
 \label{eq:derivedbalance}
\end{equation}
Smoothness of $u^*$ permits this choice of $\ell$. $\invN\ge1$ gives the other limit.
\end{proof}

\begin{proof}[Proof of Theorem~\ref{thm:discrete}]
Choose $\rho>0$ so that $C_Z\rho<(q-\gamma)/2$ and the closed $X$-ball
$\mathcal B_*$ of radius $2\rho$ centered at $u^\ast$ is contained in $\mathcal U$.
By \eqref{eq:approxinverse}, for large $N$,
\[
 C_Z\invN\epsN<(q-\gamma)/2,\quad
 C_Y\epsN<(1-q)r_N,\quad \epsN<\rho.
\]
For $u\in\mathcal B_N$, Lemma~\ref{lem:discretebounds} gives
\begin{align}
 \norm{u-u^\ast}_X&\le\epsN+\rho/\invN<2\rho,
 \label{eq:distenceBN} \\
 \norm{DS_N(u)}_{\dagger,N}
 &\le\gamma+C_Z\invN\epsN+C_Z\rho<q.
 \label{eq:nolinearDSN}
\end{align}
It holds $\mathcal B_N \subset \mathcal{B}_* \subset \mathcal{U}$ by \eqref{eq:distenceBN}. 
Since $S_N$ is $C^1$, integration along segments with \eqref{eq:nolinearDSN} 
proves it is a $q$-contraction. By \eqref{eq:defect}, we have
\[
 \norm{S_N(u_N)-u_N^\ast}_{\dagger,N}
 \le q\norm{u_N-u_N^\ast}_{\dagger,N}+C_Y\epsN<r_N.
\]
Thus $\mathcal B_N$ is invariant under $S_N$.
Banach's theorem yields the unique fixed point $u_N^\infty$, with
$\norm{u_N^\infty-u_N^\ast}_{\dagger,N}\le C_Y\epsN/(1-q)$ and
\[
 \norm{u_N^\infty-u^\ast}_X
 \le\left(1+\frac{C_Y}{1-q}\invN\right)\epsN.
\]
This proves \eqref{eq:discreteerrorX}. 
The contraction of $S_N$ in $V_N$ implies
\begin{equation}
 \norm{u_N^n-u_N^\infty}_\sigma
 \le C_aq^n\norm{u_N^0-u_N^\infty}_\sigma.
 \label{eq:contractLocalDiscrete}
\end{equation}
The split $u_N^n-u^\ast=(u_N^n-u_N^\infty)+(u_N^\infty-u^\ast)$
and the inverse estimate separate iteration and spatial error in \eqref{eq:totalerror}.
Finally, $F,\Psi,m$ vanish at $u^\ast$ and have locally bounded
derivatives. Their values at $u_N^\infty$ are therefore bounded by
$C\norm{u_N^\infty-u^\ast}_X$, which gives \eqref{eq:discreteresidual}.
\end{proof}
\begin{corollary}[Superalgebraic spatial approximation]\label{cor:cap}
Under the hypotheses of Theorem~\ref{thm:discrete}, for every $p>0$,
\begin{equation}
 \norm{u_N^\infty-u^\ast}_{C^{2,\alpha}(\overline\Omega)}
 \le C_p(N+1)^{-p},
 \label{eq:superalgebraicfixedpoint}
\end{equation}
and the residuals in \eqref{eq:discreteresidual} obey the same rate.
\end{corollary}
\begin{proof}
The preceding estimates give 
$\invN^2\epsN\le C_\ell(N+1)^{39/2-\ell}\norm{\widehat u_*}_{H^\ell}$, 
which tends to zero for $\ell>39/2$.
Substitution into \eqref{eq:discreteerrorX} and \eqref{eq:discreteresidual} 
gives the exponent $27/2-\ell$. 
Smoothness permits any integer $\ell>\max\{39/2,p+27/2\}$.
\end{proof}

For initial states $u_N^0\in\mathcal B_N$,
\begin{equation}
 \norm{u_N^n-u^\ast}_X
 \le\frac{C\rho}{\invN}q^n+C_p(N+1)^{-p}.
 \label{eq:jointconvergence}
\end{equation}
In particular $u_N^{n(N)}\to u^\ast$ in $X$ whenever $n(N)\to\infty$.
The choice $n(N)\ge p\log(N+1)/|\log q|$ gives an
$O((N+1)^{-p})$ total error.

\label{page:discrete-end}
\begingroup
\setlength{\textfloatsep}{9pt plus 2pt minus 2pt}
\setlength{\floatsep}{8pt plus 2pt minus 2pt}
\setlength{\intextsep}{8pt plus 2pt minus 2pt}
\setlength{\abovecaptionskip}{5pt}
\setlength{\belowcaptionskip}{2pt}
\renewcommand{\topfraction}{0.9}
\renewcommand{\bottomfraction}{0.8}
\renewcommand{\textfraction}{0.08}
\renewcommand{\floatpagefraction}{0.8}
\makeatletter
\let\NumericalMakeCaption\@makecaption
\long\def\@makecaption#1#2{%
  \begingroup\small\NumericalMakeCaption{#1}{#2}\endgroup}
\makeatother
\section{Numerical experiments}
\label{sec:numerics}
The section begins with point-source reflection on a spherical cap, using smooth
reference solutions to examine the convergence predicted by the theory.
Next, we turn to prescribed target images to test the method beyond
this local setting. Finally, parallel-beam refraction on a square
illustrates how the same iteration can be adapted to a different
optical model and domain.
All experiments run on an AMD Ryzen 9 9955HX laptop, with
an NVIDIA GeForce RTX 5070 Ti Laptop GPU.
\subsection{Point-source reflector construction}
\label{sec:reflector_experiments}
The reflector is $\Gamma_u=\{e^{-u(x)}x:x\in\Omega\}$ on the cap
$\Omega=\{x\in\mathbb S^2:x\cdot e_3>\cos\theta_0\}$,
$\theta_0=\pi/4$. $x$ is the incident direction and $\rho=e^{-u}$
the radius. With $p=\nabla_{\mathbb S^2}u$, the reflected direction is
\begin{equation}
 s_u=\frac{(|p|^2-1)x-2p}{1+|p|^2}.
 \label{eq:reflector_direction}
\end{equation}
In the far field, $T_u=s_u$ and the generating function is
$\Gcal(x,y,z)=\log(1-x\cdot y)-z$, giving $d_u=0$.
For a near-field spherical receiver of radius $R$, the generating
function is
\[
 \Gcal_R(x,y,\ell)
 =\log(\ell-Rx\cdot y)
  -\log\frac{\ell^2-R^2}{2},
 \quad \ell>R,
\]
where $y$ is the unit direction of the receiving point and $\ell$
is the total optical path. Set $c_u=x\cdot s_u$. At the contact point,
$\ell=e^{-u}+L_u$, with $L_u=-e^{-u}c_u+ \sqrt{R^2-e^{-2u}(1-c_u^2)}$ 
the distance from the reflector
to the receiver. This gives
$d_u=(e^{-u}/L_u)\nabla u$.

\subsubsection{Smooth reference solutions}
\label{sec:reflector_smooth}
In disk coordinates $D=\{(\xi,\eta):\xi^2+\eta^2\le1\}$, 
prescribe the nonpolynomial reference potentials
\begin{align}
 \widehat u^{\ast,\rm near}&=0.02(1-\xi^2-\eta^2)^2\sin(2\xi+\eta),\label{eq:analyticreference}\\
 \widehat u^{\ast,\rm far}&=0.01\sin(2\xi+\eta)+0.005(\xi^2-\eta^2).
 \label{eq:farreference}
\end{align}
Set $g(y)=1+0.2y_1+0.1y_2$ and manufacture
$f=g(T_*)\lvert\det(\nabla T_*)\rvert$.
The near-field receiver is the opposite cap of a radius-$3$ sphere, where 
$y$ is its unit direction, $g$ is per solid angle, and
$\sigma=1$, $a_*=g(-x)$. The far-field target
$\Lambda=T_*(\Omega)$ is noncircular, with $\sigma=0$, $a_*=f/\chi$.
Both tests use $\tau=0.2$ and independently start at every degree from 
the mean-zero quadratic seed
\begin{equation}
 \widehat u_N^0(\xi,\eta)=\widehat u_{\rm seed}(\xi,\eta)
 =0.03\xi+0.02\eta+0.01(\xi^2-\eta^2).
 \label{eq:testinitial}
\end{equation}
At $N=12$, the initial
relative potential $L^2$ errors are $308.2\%$ and $143.6\%$ (near/far).
Here $E_2=\norm{u_N-u^\ast}_{L^2(\Omega)}$ is the absolute potential error and 
potential RMS is $E_2/|\Omega|^{1/2}$. Gradient and Hessian errors below
are also absolute area-weighted $L^2$ norms.
The near-field $E_2$ reaches $1.30\times10^{-8}$ at step 35 as shown in 
Figure~\ref{fig:smooth_states}.
\begin{figure}[!htbp]
\centering\includegraphics[width=\textwidth]{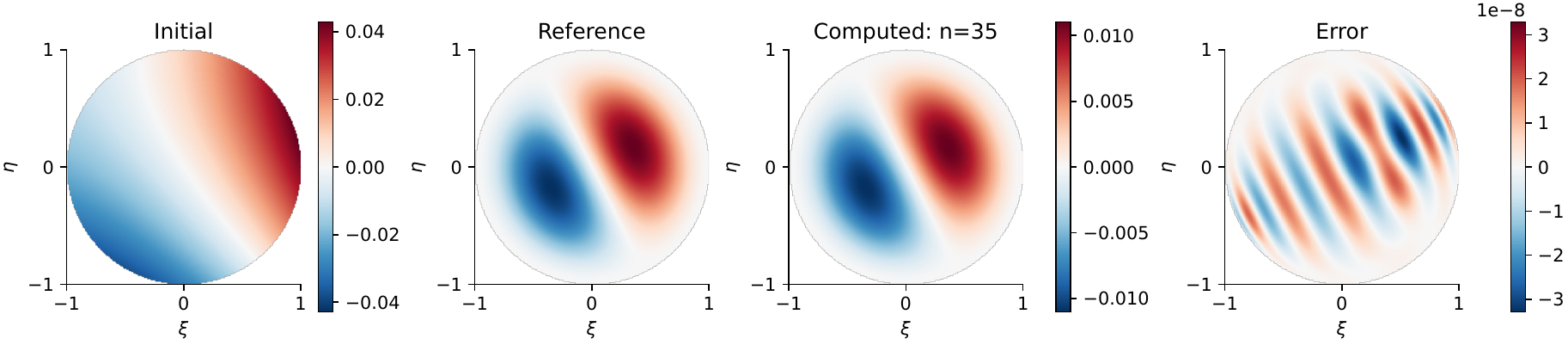}
\caption{Near-field potentials at $N=12$. The panels show the
initial potential, reference solution, stopping iterate
$u_N^{35}$, and signed pointwise error.}
\label{fig:smooth_states}
\end{figure}

Polar quadrature uses $128\times512$ (radial $\times$ angular) 
nodes for basis means and initial coefficients, and $48\times144$ for assembly. 
The interior residual $e_F$ and boundary residual $e_\partial$ 
used for stopping are defined as in Algorithm~\ref{alg:fixedmatrix}. Define
\begin{equation}
 S_n=\left(e_F(u_N^n)^2+\frac{e_\partial(u_N^n)^2}{(\pi/4)^2}\right)^{1/2},\quad
 d_n=\frac{S_{n-1}-S_n}{S_{n-1}}.
 \label{eq:plateau_stop}
\end{equation}
The iteration stops if an evaluation fails or $S_n\le10^{-10}$.
From $n=10$, it also stops for $d_n<0$ (residual increase) or
$0\le d_n<10^{-3}$, with a 600-update budget. 
The nonpolynomial tests also stop on loss of sampled $J$ or $W$ positivity.
Reference errors are diagnostics, evaluated independently on a
$113\times347$ grid and different boundary points.
\begin{figure}[!htbp]
\centering\includegraphics[width=0.88\textwidth]
{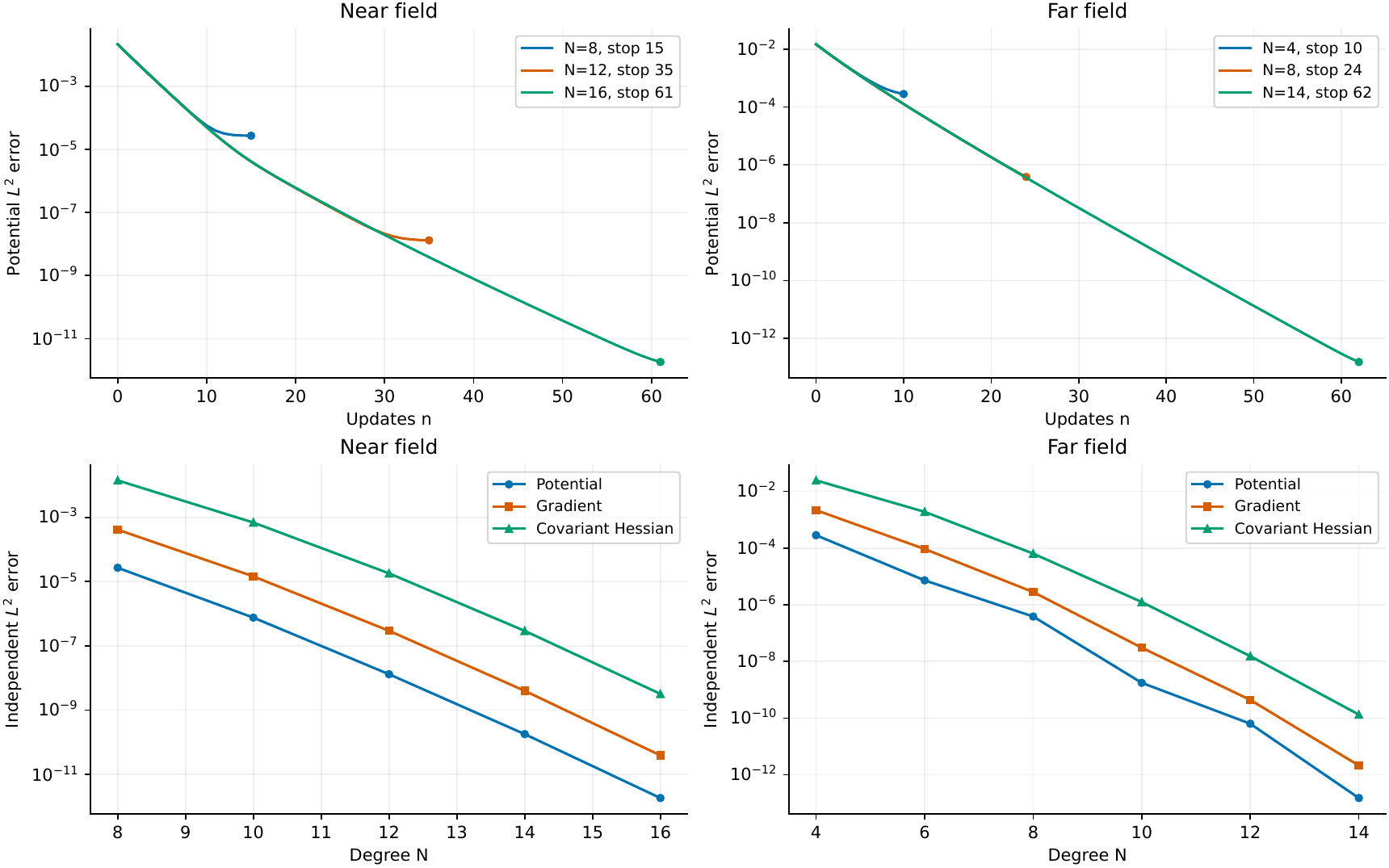}
\caption{Smooth reflector convergence. Top: potential $L^2$ errors through
the stopping states (circles). Bottom: independent potential, gradient,
and covariant-Hessian errors versus degree at the stopping states.}
\label{fig:smooth_convergence}
\end{figure}

Let $J_{T_u}^{\rm sgn}$ denote the analytic signed expression for the Jacobian,
which equals $\lvert\det(\nabla T_u)\rvert$ on the prescribed regular branch. 
All primary smooth stopping states have sampled 
$J_{T_u}^{\rm sgn}>0$ and $W\succ0$. 
The upper panels of Figure~\ref{fig:smooth_convergence} 
show approximately geometric decay of the potential error 
over the main iteration regime. 
The lower panels also show approximately geometric decay of the potential,
gradient, and covariant-Hessian errors as $N$ increases.
The potential error decreases from $2.69\times10^{-5}$ to 
$1.82\times10^{-12}$ for near-field $N=8-16$ 
and from $2.83\times10^{-4}$ to $1.52\times10^{-13}$ for far-field, $N=4-14$. 
Table~\ref{tab:refinement} reports the 
corresponding stopping steps and residuals.
Table~\ref{tab:smooth_sigma} shows that 
changing $\sigma$ affects the iteration speed and stopping error.
For both models, $\sigma=5$ requires more
updates and gives larger errors than $\sigma=0$ or $1$.
The near-field test also converges at $\sigma=0$ outside the
theorem's sufficient parameter range.
\begin{table}[!htbp]
\centering
\begin{minipage}[t]{0.55\textwidth}
\vspace{0pt}\centering\small
\setlength{\tabcolsep}{4pt}
\caption{Stopping steps $n_s$ and independent residuals.}
\label{tab:refinement}
\begin{tabular}{lrrrrr}
\toprule
Model & $N$ & $d_N$ & $n_s$ & $e_F$ & $e_\partial$\\
\midrule
\InputTableRows{data_adjacent/smooth_validation.tex}
\bottomrule
\end{tabular}
\end{minipage}\hfill
\begin{minipage}[t]{0.42\textwidth}
\vspace{0pt}\centering\small
\setlength{\tabcolsep}{4pt}
\caption{Shift comparison at $N=12$: 
stopping steps and potential $L^2$ errors $E_2$, scaled as indicated.}
\label{tab:smooth_sigma}
\begin{tabular}{rrrrr}
\toprule
 & \multicolumn{2}{c}{Near field} & \multicolumn{2}{c}{Far field}\\
$\sigma$ & $n_s$ & $10^8E_2$ & $n_s$ & $10^{11}E_2$\\
\midrule
\InputTableRows{data_adjacent/sigma.tex}
\bottomrule
\end{tabular}
\end{minipage}
\end{table}

Figure~\ref{fig:boundary}(a) shows that the iteration corrects the
mapped boundary, reducing its RMS error to $1.27\times10^{-11}$
radians in 62 updates at $N=14$.
Panel (b) is a separate boundary-load ablation with the polynomial reference
$\widehat u^\ast=0.02\xi(1-\xi^2-\eta^2)^2$ of degree five, exactly representable at $N=8$.
It uses the same quadratic seed, $R=3$, $g(y)=1+0.2y_1+0.1y_2$,
manufactured $f=g(T_*)J_*$, and matched $a_*=g(-x)$.
Set $\sigma=0$, $\tau=0.2$, $32\times96$ quadrature and a 100-update budget. 
The two runs differ only by retaining or removing the boundary load.
With the boundary term, $S_n\le10^{-10}$ at step 40.
Without it, the iteration fails to converge.
\begin{figure}[!htbp]
\centering\includegraphics[width=0.88\textwidth]{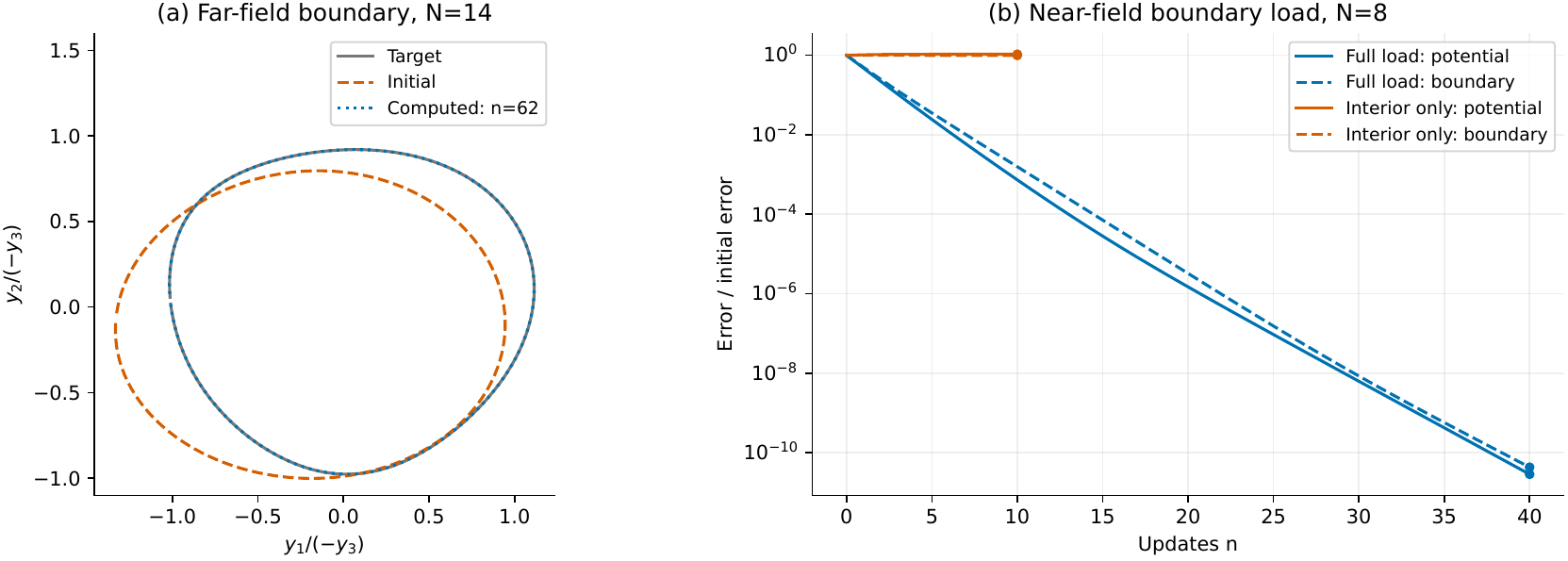}
\caption{
(a) Noncircular far-field target and mapped initial and stopping boundaries, $N=14$. 
(b) Potential RMS (solid) and boundary RMS (dashed) for the separate
polynomial near-field test at $N=8$,
normalized by their initial values, with (blue) and without
(orange) the boundary correction. Circles mark stopping states.}
\label{fig:boundary}
\end{figure}

\subsubsection{Discontinuous far-field targets}
\label{sec:reflector_far_images}
The $128\times128$ disk targets in Figure~\ref{fig:raw_degree} are mapped
to a receiving spherical cap by central projection. 
The image intensities are converted from density per unit planar
area to density per unit solid angle. The PDE target density is normalized
to the total source mass $M=|\Omega|=2\pi(1-\cos\theta_0)$ for $f=1$,
so $\int_\Lambda g\dd V_\Lambda=M$.
Use $f=1$, $N=128$, $256\times1024$ quadrature, 1024 boundary points,
$u^0=0$, $\sigma=0$, $\tau=0.2$, and frozen $a=a_0=1$ calibrated at the
initial state. 
Rule~\eqref{eq:plateau_stop} uses a 500-update budget.
Weak residuals are only recorded. 
The runs stop at 26/20 updates for the two cases. 
For independent validation, we trace $2^{22}$ equal-energy 
rays and record their hits.
Let $H_j$ be the fraction of emitted energy received by pixel $j$,
and $p_j$ be its prescribed target energy fraction, normalized separately
by $M$ so that $\sum_j p_j=1$.
Set $H=(H_j)$ and $p=(p_j)$. Define the captured energy fraction
and relative shape error by
\begin{equation}
 c=\sum_jH_j,\quad E_{\rm shape}=\frac{\norm{H/c-p}_{\ell^2}}{\norm p_{\ell^2}}.
 \label{eq:shape_metric}
\end{equation}
\begin{figure}[!htbp]
\centering\includegraphics[width=\textwidth]{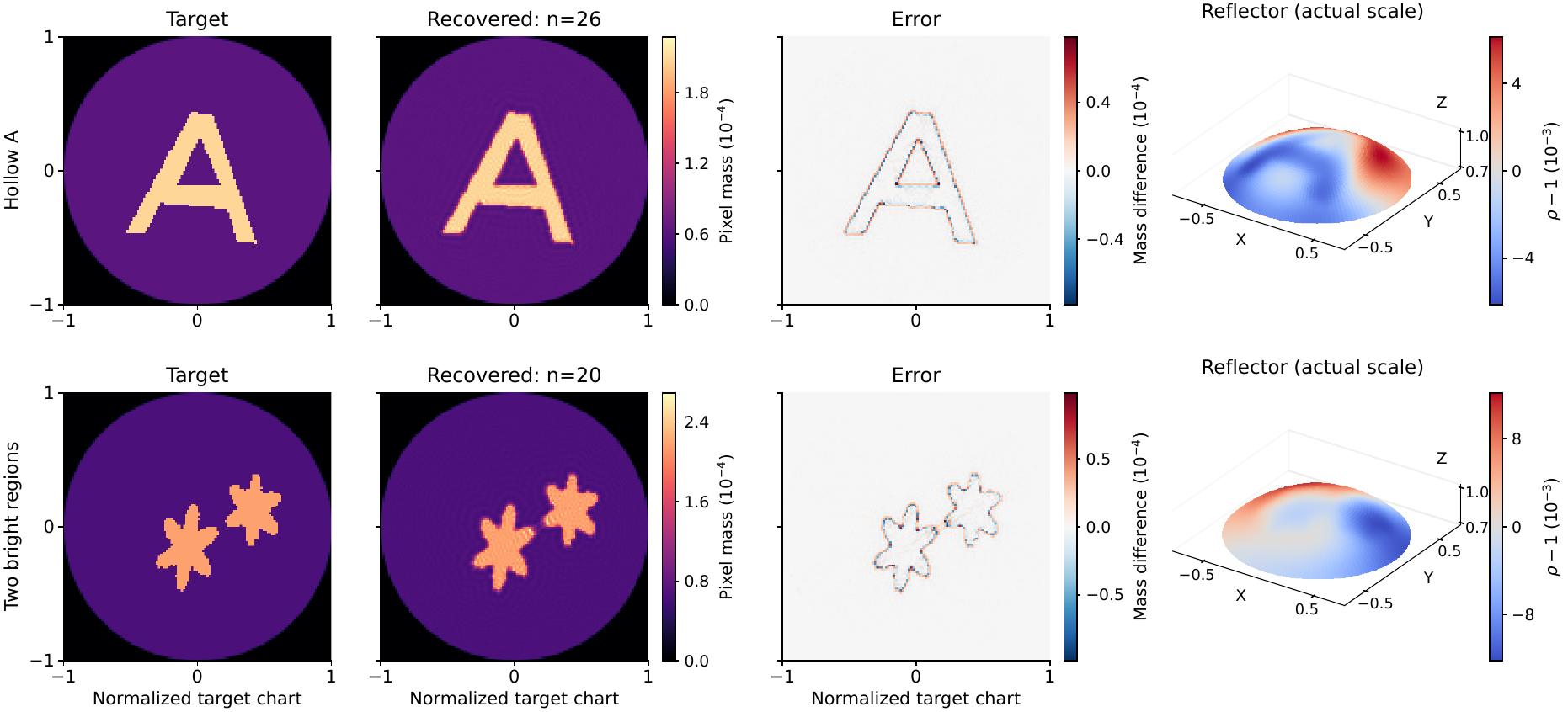}
\caption{Far-field targets, ray-traced images, errors $H-p$, and surfaces.
Surface coordinate units are equal, with color $\rho-1$.}
\label{fig:raw_degree}\label{fig:raw_two}
\end{figure}

\subsubsection{Complex near-field targets}
\label{sec:reflector_near_images}
Near-field rays meet a radius-$3$ disk at $z=-H$.
The target density is measured per unit physical area and normalized to
the total source mass $M$. Writing $\rho_u=e^{-u}$, the reflected ray
$\rho_u x+Ls_u$ intersects the receiver at
\[
 L_u=\frac{-H-\rho_u x_3}{(s_u)_3}>0,\quad
 T_u=\rho_u x_{1,2}+L_u(s_u)_{1,2}.
\]
Here $x_{1,2}=(x_1,x_2)$, and $(s_u)_{1,2}$ denotes the first two
components of the reflected direction. The boundary distance on the
receiver plane is $\psi(y)=|y|-3$.
Set $\sigma=1$, $\tau=0.02$, the practical frozen weight $a_0=1/(2H)$,
and $u^0=0$. The stopping rule replaces the angular scale $\pi/4$ 
by the receiver radius $3$. 
Figure~\ref{fig:raw_near} shows the letter A at $H=12$
stops for insufficient decrease in 177 updates, and the M--T--rectangle
target at $H=6$ stops for residual increase after 115.

Table~\ref{tab:imaging} compares CPU and GPU runtimes in double
precision, using one CPU thread. 
The GPU offers a significant speed advantage for 
computationally demanding tasks.
The table also reports the relative shape error
$E_{\rm shape}$, captured energy fraction $c$, and source-area
fraction $\eta_J$ where $J_{T_u}^{\rm sgn}\le0$, all as percentages.
\begingroup
\setlength{\intextsep}{4pt}
\setlength{\abovecaptionskip}{4pt}
\setlength{\belowcaptionskip}{0pt}
\begin{table}[!htbp]
\centering\small
\caption{Reflector stopping states and paired CPU/GPU times (seconds).}
\label{tab:imaging}\label{tab:raw_degree}\label{tab:raw_near}
\setlength{\tabcolsep}{4pt}
\begin{tabular}{lrrrrrrr}
\toprule
Case & Steps & $E_{\rm shape}/\%$ & $c/\%$ & $e_F$ & $\eta_J/\%$ & CPU & GPU\\
\midrule
\InputTableRows{data_adjacent/imaging_summary.tex}
\bottomrule
\end{tabular}
\end{table}

\begin{figure}[!htbp]
\centering\includegraphics[width=\textwidth]{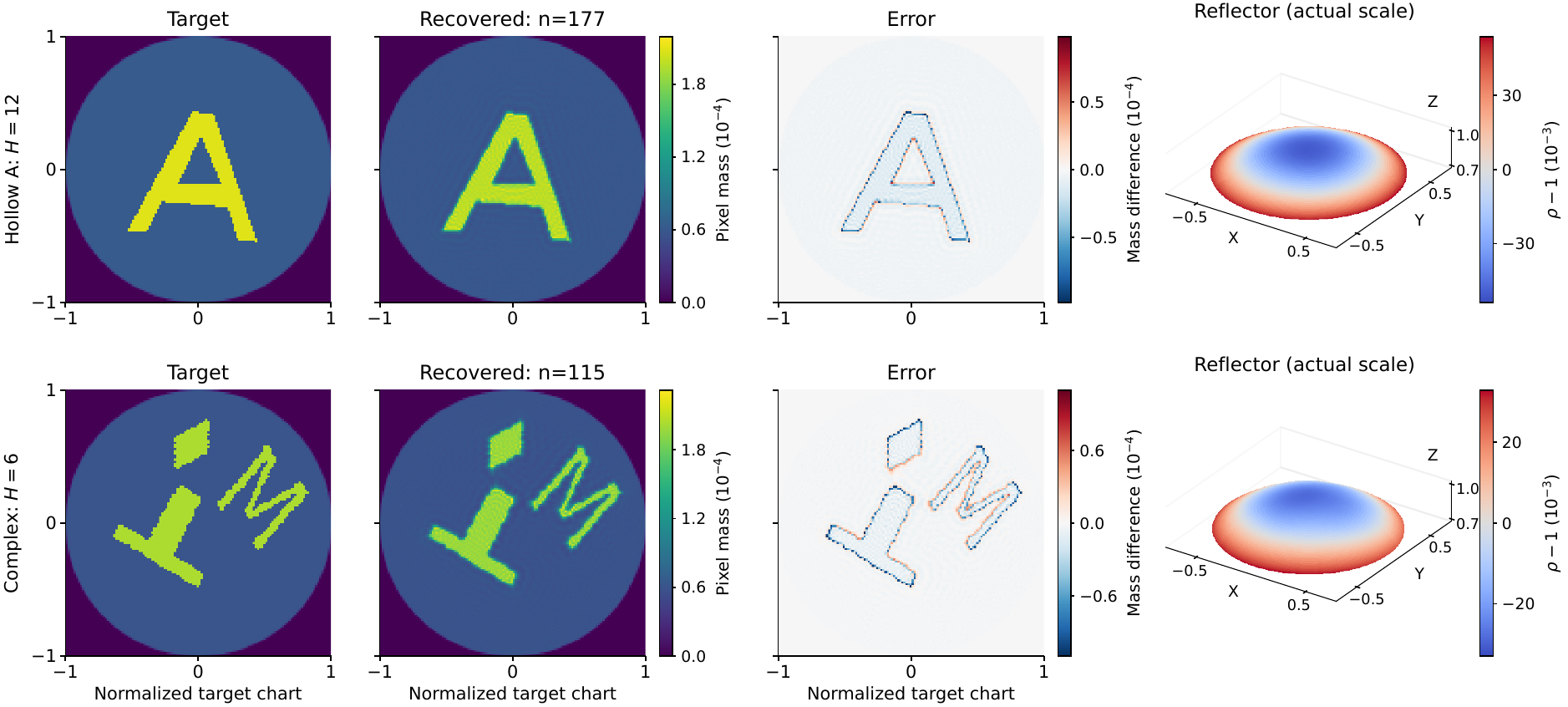}
\caption{Near-field targets, images, errors, and surfaces. 
Letter A stops at 177 updates for $H=12$, and the complex
target stops at 115 updates for $H=6$.}
\label{fig:raw_near}
\end{figure}
\endgroup

The computed reflector surfaces reproduce the main target features 
in all four image-based tests, 
with errors concentrated near image boundaries. 
These results demonstrate useful imaging performance for 
discontinuous target densities. 
They also suggest that the flux-matched Poisson iteration can 
be applied to a broader class of optical design problems 
beyond the smooth-target setting considered in the convergence analysis.

\Needspace{5\baselineskip}
\subsection{Parallel-beam refractor construction}
\label{sec:square_refractor}
We next apply the iteration on a square to design a refracting surface
$\Gamma_u=\{(x,h(x)):x\in\Omega\}$, with $h=1+u$ and $\mean u=0$. 
All lengths are scaled by $10\,\mathrm{mm}$.
Parallel rays travel along $+e_3$ from glass of refractive index
$\eta=1.5$ into air of refractive index $1$, then reach the plane $z=H=4$.
The source and target domains are $\Omega=\Lambda=(-1,1)^2$. 
The source density is $f=1$, and both densities are per unit planar area.

For a receiving point $y$, the optical path is
$\ell=\eta h+\sqrt{|x-y|^2+(H-h)^2}$.
With $\gamma=\eta^2-1$, solving for $h$ and using $u=h-1$ gives the generator
\begin{equation}
 \Gcal(x,y,\ell)=
 \frac{\eta\ell-H-\sqrt{(\ell-\eta H)^2+\gamma|x-y|^2}}{\gamma}-1.
 \label{eq:refractor_generator}
\end{equation}
Set $p=\nabla u$, $q=\sqrt{1-\gamma|p|^2}$ and
$\beta=\frac{(H-h)\gamma}{1+\eta q}$.
The ray map and ellipticity matrix are
\begin{equation}
 T_u=x+\beta p,\quad
 W_u=D^2u-\Gcal_{xx}=D^2u+\frac{I-\gamma p\otimes p}{\beta}.
 \label{eq:refractor_map}
\end{equation}
We solve the energy-conservation equation $F(u)=1-g(T_u)\det DT_u=0$
on the branch where \mbox{$0<h<H$}, $q>0$, $W_u\succ0$.
The boundary condition $\partial_\nu u=0$ keeps each side mapped to itself.

Use the tensor cosine modes
$\phi_{ij}(x)=\cos \bigl(i\pi(x_1+1)/2\bigr)
\cos \bigl(j\pi(x_2+1)/2\bigr)$.
Omitting the constant mode gives the zero-mean space
\begin{equation}
 V_N^\square=\operatorname{span}\{\phi_{ij}:0\le i,j\le N,\ (i,j)\ne(0,0)\},
 \quad \dim V_N^\square=(N+1)^2-1.
 \label{eq:squarebasis}
\end{equation}
The cosine basis satisfies the homogeneous Neumann condition,
so no boundary correction is needed.
On a uniform $(M+1)^2$ grid, $P_N^M$ retains modes
$0\le i,j\le N$ using an FFT-based cosine transform.
Derivatives are evaluated spectrally, and the Poisson solve
divides each nonconstant coefficient by
$\sigma+\pi^2(i^2+j^2)/4$.

We use two targets, both normalized to mean one.
The Gaussian target has amplitudes $(1,0.7)$, centers
$(-0.38,0.28)$ and $(0.35,-0.30)$, standard deviations
$(0.18,0.25)$ and $(0.25,0.14)$, and background $0.25$,
with even extension across the boundaries.
For the ZJU image, we assign densities
$1$ to the letters and $0.3$ to the background, then apply Gaussian
smoothing with standard deviation one pixel.

Set $N=128$, $M=256$, and $\sigma=0$. Starting from $u=0$,
solve $-\Delta v=P_N^MF-\mean{P_N^MF}$ with homogeneous Neumann
conditions and update $u$ by $u-0.25v$.
We monitor the full residual $R_M=\|F\|_M$ and record the projected residual
$r_N=\|P_N^MF-\mean{P_N^MF}\|_N$.
Here $\|\cdot\|_K$ is the normalized trapezoidal RMS norm
on a $(K+1)^2$ grid.

The stopping condition $R_M\le10^{-10}$ is checked from the initial state.
From step 10, it also stops if $R_M$ increases
or decreases by less than $0.5\%$.
The budget is 500 updates. 
\begin{table}[!htbp]
\centering\small
\setlength{\tabcolsep}{5pt}
\begin{tabular}{lrrrrr}
\toprule
Target & Steps & $R_M$ & $r_N$ & Independent RMS & $100E_{\rm shape}$\\
\midrule
\InputTableRows{data_adjacent/refractor_results.tex}
\bottomrule
\end{tabular}
\caption{Refractor stopping residuals and image error \eqref{eq:shape_metric}.}
\label{tab:refractor_results}
\end{table}

\begin{figure}[!htbp]
\centering\includegraphics[width=\textwidth]{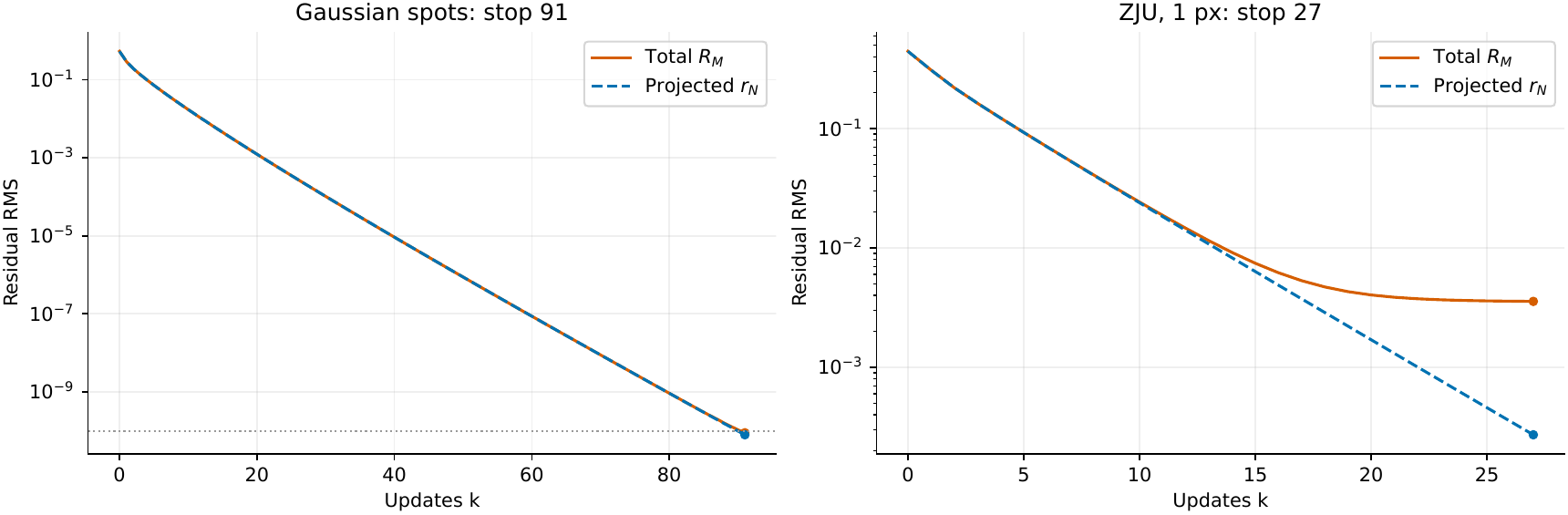}
\caption{Refractor total and projected residuals through the stopping states.}
\label{fig:refractor_residuals}
\end{figure}

\begin{figure}[!htbp]
\centering\includegraphics[width=\textwidth]{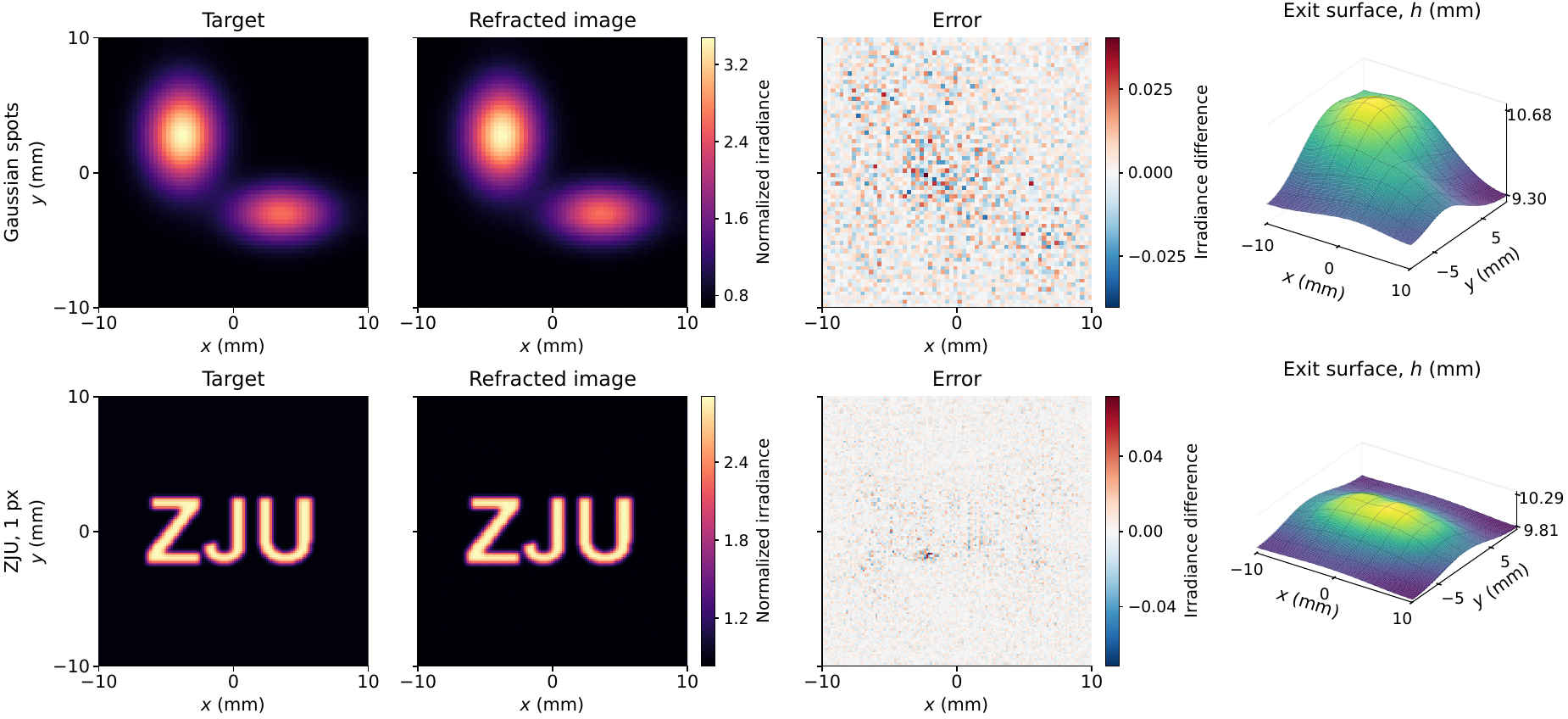}
\caption{Refractor targets, vector-Snell images, signed errors, and exit
heights. }
\label{fig:refractor_images}
\end{figure}

Table~\ref{tab:refractor_results} shows that the Gaussian case reaches
a residual below $10^{-10}$ in 91 updates.
The nearly linear semilog curve in Figure~\ref{fig:refractor_residuals}
is consistent with geometric decay in this experiment, which lies outside
the present convergence theory.
For ZJU, the full residual levels off and triggers stopping after
27 updates, while the projected residual continues to decrease.
Both image errors are below $0.71\%$.
Figure~\ref{fig:refractor_images} shows that the constructed surfaces
reproduce both target patterns, demonstrating the method's
effectiveness for these square-domain refractor designs.

\FloatBarrier
\endgroup

\FloatBarrier
\section{Conclusion}
\label{sec:conclusion}

In this paper, we proposed a flux-matched Poisson 
iteration for solving generated
Jacobian equations with second boundary conditions.
Each correction solves a fixed shifted Poisson problem with a boundary
correction. Near a smooth regular solution, we proved local
geometric convergence under suitable assumptions.
For the Fourier--Jacobi discretization on spherical caps, we
established discrete local convergence and error bounds that
separate iteration and spatial errors.

The smooth numerical tests also support these results, while the design
experiments illustrate empirical optical performance.
Smooth tests show geometric convergence and rapid decay of spatial
errors. Image-based tests reproduce the main target features,
including cases beyond the scope of the theory. The signed reflector
stopping states are not certified regular GJE solutions.
The square-domain tests also show that the method can be used
with other domains and bases.

The image-based tests suggest that the method can work beyond
the smooth setting studied here. This motivates further study
of discrete targets and less regular densities, with convergence
analysis in weaker function spaces. For practical use, it would
also be useful to choose boundary weights without a known reference
solution and to account for numerical quadrature errors.

\label{page:body-end}

\end{document}